\input ./style/arxiv-general.cfg
\documentclass[MSNbibl,number,citesort,seceqn,dvips]{arxbj}
\makeatletter
   \@ifpackageloaded{graphicx}{}{\usepackage{graphicx}}
\makeatother
\usepackage{mathrsfs}
\aid{0}
\volume{21}
\issue{3}
\pubyear{2015}
\firstpage{1670}
\lastpage{1696}
\doi{10.3150/14-BEJ617} % kopijuoti is 'New paper accepted'
\makeatletter
\newcommand{\rrvert}{\vert}
\newcommand{\llvert}{\vert}
\def\id{\mathbf{id}}
\def\idn{\mathbf{id}_n }
\def\idtwo{\mathbf{id}_2 }
\def\idij{\mathbf{id}_{\{i,j\}} }
\def\wireN{\mathcal{W}_{\infty}}
\def\wirem{\mathcal{W}_m}
\def\Rmn{\mathbf{R}_{m,n}}

\def\symmetricN{\mathscr{S}_{\mathbb{N}}}

\def\graphsm{\mathcal{G}_m}

\def\graphsN{\mathcal{G}_{\infty}}
\def\graphlimits{\mathcal{D}^*}
\def\rewiringlimits{\mathcal{V}^*}
\def\ind{\operatorname{ind}}

\def\equalinlaw{=_{\mathcal{L}}}
\newcommand{\eqref}[1]{(\ref{#1})}
\newtheorem{theorem}{Theorem}[section]
\newremark{rmk}{Remark}[section]
\newproclaim{example}{Example}[section]
\newtheorem{lemma}{Lemma}[section]
\newtheorem{prop}{Proposition}[section]
\newtheorem{cor}{Corollary}[section]
\newproclaim{defn}{Definition}[section]
\makeatother

\begin{document}
\begin{frontmatter}

\title{Time-varying network models}
\runtitle{Time-varying network models}

\begin{aug}
%%%% inicialai - be tarpu
\author[1]{\inits{H.}\fnms{Harry} \snm{Crane}\corref{}\ead[label=e1]{hcrane@stat.rutgers.edu}}% \and
%\author{\inits{}\fnms{}~\snm{}\thanksref{}\ead[label=e2]{}}
%\author{\inits{}\fnms{}~\snm{}}
%%\runauthor{} %% auto
%\dedicated{}
\address[1]{Department of Statistics \& Biostatistics, Rutgers University,
110 Frelinghuysen Road,
Hill Center, Room~461,
Piscataway, NJ 08854, USA. \printead{e1}}
%\address[]{}
\end{aug}

% HISTORY:
\received{\smonth{6} \syear{2013}}
\revised{\smonth{1} \syear{2014}}

% ABSTRACT
%
\begin{abstract}
We introduce the {\em exchangeable rewiring process} for modeling
time-varying networks.
The process fulfills fundamental mathematical and statistical
properties and can be easily constructed from the novel operation of
{\em random rewiring}.
We derive basic properties of the model, including consistency under
subsampling, exchangeability, and the Feller property.
A reversible sub-family related to the Erd\H{o}s--R\'enyi model arises
as a special case.
\end{abstract}

% KEYWORDS
% visi is mazosios raides ir pagal abecele
%

\end{frontmatter}

%s1 #&#
\section{Introduction}\label{section:introduction}

A recent influx of academic monographs \cite
{ChungLubook,DorogovtsevMendes2003,DurrettRandomGraphs,Kolaczykbook,LovaszBook,vanderHofstadBook}
and popular books \cite{Linked,ClarkeCyber,SixDegrees} manifests a
keen cultural and scientific interest in complex networks, which appeal
to both applied and theoretical problems in national defense,
sociology, epidemiology, computer science, statistics, and mathematics.
The Erd\H{o}s--R\'enyi random graph \cite{ErdosRenyi1959,ErdosRenyi1960} remains the most widely studied network model.
Its simple dynamics endow it with remarkable mathematical properties,
but this simplicity overpowers any ability to replicate realistic structure.
Many other network models have been inspired by empirical observations.
Chief among these is the {\em scale-free} phenomenon, which has
garnered attention since the initial observation of power law behavior
for Internet statistics \cite{Faloutsos1999}.
Celebrated is Barab\'asi and Albert's preferential attachment model
\cite{BarabasiAlbert1999}, whose dynamics are tied to the {\em rich
get richer} or {\em Matthew effect}.\footnote{``For to everyone who
has will more be given, and he will have an abundance. But from the one
who has not, even what he has will be taken away.'' (Matthew 25:29,
{\em The Bible}, English Standard Version, 2001.)}
Citing overlooked attributes of network sampling schemes, other authors
\cite{LiAldersonDoyleWillingerGraphs,WillingerAldersonDoyle2009} have
questioned the power law's apparent ubiquity.
Otherwise, Watts and Strogatz \cite{WattsStrogatz1998} proposed a
model that replicates Milgram's {\em small-world} phenomenon \cite
{Milgram1967}, the vernacular notion of {\em six degrees of separation}
in social networks.

Networks arising in many practical settings are dynamic, they change
with time.
Consider a population $\{u_1,u_2,\ldots\}$ of individuals.
For each $t\geq0$, let $G_{ij}(t)$ indicate a social relationship
between $u_i$ and $u_j$ and let $G_t:=(G_{ij}(t))_{i,j\geq1}$ comprise
the indicators for the whole population at time $t$.
For example, $G_{ij}(t)$ can indicate whether $u_i$ and $u_j$ are
co-workers, friends, or family, have communicated by phone, email, or
telegraph within the last week, month, or year, or subscribe to the
same religious, political, or philosophical ideology.
Within the narrow scope of social networks, the potential meanings of
$G_{ij}(t)$ seem endless; expanding to other disciplines, the possible
interpretations grow.
In sociology, $\{G_t\}_{t\geq0}$ records changes of social
relationships in a population; in other fields, the network dynamics
reflect different phenomena and, therefore, can exhibit vastly
different behaviors.
In each case, $\{G_t\}_{t\geq0}$ is a time-varying network.

Time-varying network models have been proposed previously in the
applied statistics literature.
The {\em Temporal Exponential Random Graph Model} (TERGM) in \cite
{HannekeFuXing2010} incorporates temporal dependence into the {\em
Exponential Random Graph Model} (ERGM).
The authors highlight select properties of the TERGM, but consistency
under subsampling is not among them.
From the connection between sampling consistency and lack of
interference, it is no surprise that the Exponential Random Graph Model
is sampling consistent only under a choking restriction on its
sufficient statistics \cite{RinaldoShalizi2013}.
McCullagh \cite{McCullagh2002} argues unequivocally the importance of
consistency for statistical models.

Presently, no network model both meets these logical requirements and
reflects empirical observations.
In this paper, rather than focus on a particular application, we
discuss network modeling from first principles.
We model time-varying networks by stochastic processes with a few
natural invariance properties, specifically, exchangeable, consistent
Markov processes.
%The rewiring maps are vital to our construction and are introduced in
%Section~\ref{section:rewiring maps}.

The paper is organized as follows. In Section~\ref{section:modeling
preliminaries}, we discuss first principles for modeling time-varying
networks; in Section~\ref{section:informal description}, we describe
the rewiring process informally; in Section~\ref{section:rewiring
maps}, we introduce the workhorse of the paper, the rewiring maps; in
Sections~\ref{section:discrete} and \ref{section:exchangeable
rewiring maps}, we discuss a family of time-varying network models in
discrete-time; in Section~\ref{section:continuous}, we extend to
continuous-time; in Section~\ref{section:Poissonian structure}, we
show a Poisson point process construction for the rewiring process, and
we use this technique to establish the Feller property; and in Section~\ref{section:concluding remarks}, we make some concluding remarks. We
prove some technical lemmas and theorems in Section~\ref{section:proof}.

%s2 #&#
\section{Modeling preliminaries}\label{section:modeling preliminaries}

For now, we operate with the usual definition of a graph/network as a
pair $G:=(V,E)$ of vertices and edges.
We delay formalities until they are needed.

Let $\boldsymbol{\Gamma}:=\{\Gamma_t\}_{t\in T}$ be a random collection
of graphs indexed by $T$, denoting {\em time}.
We may think of $\boldsymbol{\Gamma}$ as a collection of social networks
(for the same population) that changes as a result of social forces,
for example, geographical relocation, broken relationships, new
relationships, etc., but our discussion generalizes to other applications.

In practice, we can observe only a finite sample of individuals.
Since the population size is often unknown, we assume an infinite
population so that our model only depends on known quantities.
Thus, each $\Gamma_t$ is a graph with infinitely many vertices, of
which we observe a finite sub-network $\Gamma^{[n]}_t$ with
$n=1,2,\ldots$ vertices.
Since the vertex labels play no role, we always assume sampled graphs
have vertex set $[n]:=\{1,\ldots,n\}$, where $n$ is the sample size,
and the population graph is infinite with vertex set $\mathbb{N}:=\{
1,2,\ldots\}$, the natural numbers.

The models we consider are {\em Markovian}, {\em exchangeable}, and
{\em consistent}.

%s2.1 #&#
\subsection{The Markov property}\label{section:Markov property}

The process $\boldsymbol{\Gamma}$ has the {\em Markov property} if, for
every $t>0$, its pre-$t$ and post-$t$ $\sigma$-fields are
conditionally independent given the present state $\Gamma_t$. Put
another way, the current state $\Gamma_t$ incorporates all past and
present information about the process, and so the future evolution
depends on $\sigma\langle\Gamma_s\rangle_{s\leq t}$ only through
$\Gamma_t$.

It is easy to conceive of counterarguments to this assumption: in a
social network, suppose there is no edge between individuals $i$ and
$i'$ or between $j$ and $j'$ at time $t>0$.
Then, informally,\footnote{We are implicitly ignoring the dependence
between $ii'$ and $jj'$ for the sake of illustration.} we expect the
future (marginal) evolution of edges $ii'$ and $jj'$ to be identically
distributed.
But if, in the past, $i$ and $i'$ have been frequently connected and
$j$ and $j'$ have not, we might infer that the latent relationships
among these individuals are different and, thus, their corresponding
edges should evolve differently. For instance, given their past
behavior, we might expect that $i$ and $i'$ are more likely than $j$
and $j'$ to reconnect in the future.

Despite such counterarguments, the Markov property is widely used and
works well in practice. Generalizations to the Markov property may be
appropriate for specific applications, but they run the risk of overfitting.

%s2.2 #&#
\subsection{Exchangeability}\label{section:exchangeability}
Structure and changes to structure drive our study of networks.
Vertex labels carry no substantive meaning other than to keep track of
this structure over time; thus, a suitable model is {\em exchangeable},
that is, its distributions are invariant under relabeling of the vertices.

For a model on finite networks (i.e., finitely many vertices),
exchangeability can be induced trivially by averaging uniformly over
all permutations of the vertices.
But we assume an infinite population, for which the appropriate
invariance is {\em infinite exchangeability}, the combination of
exchangeability and consistency under subsampling (Section~\ref{section:consistency}).
Unlike the finite setting, infinite exchangeability cannot be imposed
arbitrarily by averaging; it must be an inherent feature of the model.

%s2.3 #&#
\subsection{Markovian consistency under subsampling}\label
{section:consistency}
For any graph with vertex set $V$, there is a natural and obvious
restriction to an induced subgraph with vertex set $V'\subset V$ by
removing all vertices and edges that are not fully contained in $V'$.
The assumption of {\em Markovian consistency}, or simply {\em
consistency}, for a graph-valued Markov process implies that, for every
$n\in\mathbb{N}$, the restriction $\boldsymbol{\Gamma}^{[n]}$ of
$\boldsymbol{\Gamma}$ to the space of graphs with vertex set $[n]$ is,
itself, a Markov process. Note that this property does not follow
immediately from the Markov assumption for $\boldsymbol{\Gamma}$ because
the restriction operation is a many-to-one function and, in general, a
function of a Markov process need not be Markov. Also note that the
behavior of the restriction $\boldsymbol{\Gamma}^{[n]}$ can depend on
$\boldsymbol{\Gamma}$ through as much as its exchangeable $\sigma
$-field, which depends only on the ``tail'' of the process.

Markovian consistency may be unjustified in some network modeling
applications. This contrasts with other combinatorial stochastic
process models, for example, coalescent processes \cite{Kingman1982},
for which consistency is induced by an inherent {\em lack of
interference} in the underlying scientific phenomena.
Nevertheless, if we assume the network is a sample from a larger
network, then consistency permits out-of-sample statistical inference
\cite{McCullagh2002}. Without Markovian consistency in a time-varying
Markov model, the sampled process can depend on the whole (unobserved)
process, leaving little hope for meaningful inference.

%s3 #&#
\section{Rewiring processes: Informal description}\label
{section:informal description}
We can envision at least two kinds of network dynamics that correspond,
intuitively, to {\em local} and {\em global} structural changes.
Local changes involve {\em only one} edge, global changes involve a
{\em positive fraction} of edges.
We say the {\em status} of edge $ij$ is {\em on} if there is an edge
between $i$ and $j$; otherwise, we say the status is {\em off}.

A local change occurs whenever the status of exactly one edge changes,
called a {\em single-edge update}.
An easy way to generate single-edge updates is by superposition of
independent rate-1 Poisson processes. For each pair $i<j$, we let $\{
T^{ij}_k\}_{k\geq1}$ be the arrival times of a rate-1 Poisson point
process. At each arrival time, the status of the edge between $i$ and
$j$ changes (either from `off' to `on' or the reverse). Doing this
independently for each pair results in an infinite number of changes to
the network in any arbitrary time interval, but only finitely many
changes within each finite subnetwork. We call a process with this
description a {\em local-edge process}; see Section~\ref{section:local-edge}.

A global change occurs whenever the status of a positive proportion of
edges changes simultaneously.
In practice, such an event might indicate a major external disturbance
within the population, for example, spread or fear of a pandemic.
Modeling such processes in continuous-time requires more preparation
than the local-edge process.

For an example, consider generating a discrete-time Markov chain
$\boldsymbol{\Gamma}:=\{\Gamma_m\}_{m=0,1,2,\ldots}$ on the finite
space of graphs with vertex set $[n]$. At any time $m$, given $\Gamma
_m=G$, we can generate a transition to a new state $G'$ as follows.
Independently for each pair $i<j$, we flip a coin to determine whether
to put an edge between $i$ and $j$ in $G'$: if $ij$ is {\em on} in $G$,
we flip a $p_1$-coin; otherwise, we flip a $p_0$-coin. This description
results in a simple, exchangeable Markov chain on finite graphs, which
we call the {\em Erd\H{o}s--R\'enyi rewiring chain} (Section~\ref{section:ER}). More general transitions are possible, for example,
edges need not evolve independently.
We use the next Markov chain as a running example of a discrete-time
rewiring chain.

%s3.1 #&#
\subsection{A reversible Markov chain on graphs}\label{section:example}
We fix $n\in\mathbb{N}$ and regard an undirected graph $G$ with
vertex set $[n]$ as a $\{0,1\}$-valued symmetric matrix
$(G_{ij})_{1\leq i,j\leq n}$ such that $G_{ii}=0$ for all $i=1,\ldots
,n$; that is, we represent a graph by its adjacency matrix with
$G_{ij}:=\mathbf{1}\{ij\mbox{ is on}\}$.
For any pair of graphs $(G,G')$, we can compute the statistic $\mathbf
{n}:=\mathbf{n}(G,G'):=(n_{00},n_{01},n_{10},n_{11})$, where for $r,s=0,1$,
\[
n_{rs}:=\sum_{1\leq i<j\leq n}\mathbf{1}\bigl
\{G_{ij}=r\mbox{ and }G'_{ij}=s\bigr\}.
\]
For example, $n_{01}$ is the number of pairs $i,j$ for which the status
of $ij$ changes from $0$ to $1$ from $G$ to $G'$. We use $\mathbf{n}$
as a sufficient statistic to define the transition probability
\[
P_{\alpha,\beta}^{(n)}\bigl(G,G'\bigr):=\frac{\alpha^{\uparrow n_{00}}\beta
^{\uparrow n_{01}}\alpha^{\uparrow n_{11}}\beta^{\uparrow
n_{10}}}{(\alpha+\beta)^{\uparrow(n_{01}+n_{00})}(\alpha+\beta
)^{\uparrow(n_{10}+n_{11})}},
\]
where $\alpha^{\uparrow j}:=\alpha(\alpha+1)\cdots(\alpha+j-1)$
and $\alpha,\beta>0$.

The sufficient statistic $\mathbf{n}$ is invariant under joint
relabeling of the vertices of $(G,G')$ and so the transition law is
exchangeable. Furthermore, $P_{\alpha,\beta}^{(n)}$ is reversible
with respect to
\[
\varepsilon_{\alpha+\beta,\alpha+\beta}^{(n)}(G):=\frac{(\alpha
+\beta)^{\uparrow n_0}(\alpha+\beta)^{\uparrow n_1}}{(2\alpha
+2\beta)^{\uparrow n}},
\]
where $n_r:=\sum_{1\leq i<j\leq n}\mathbf{1}\{G_{ij}=r\}$, $r=0,1$.
The distribution $\varepsilon_{\alpha,\beta}^{(n)}$ arises as a
mixture of Erd\H{o}s--R\'enyi random graphs with respect to the
$\operatorname{Beta}(\alpha,\beta)$ distribution.
Furthermore, $\{P_{\alpha,\beta}^{(n)}\}_{n\in\mathbb{N}}$ is a
consistent collection of transition probabilities and, therefore,
determines a unique transition probability (and hence Markov chain) on
the space of infinite graphs with vertex set $\mathbb{N}$.

Though consistency is not immediately obvious for the above family, the
savvy reader might anticipate it: the formula for $P_{\alpha,\beta
}^{(n)}$ involves rising factorials (i.e., Gamma functions), which
also appear in other consistent combinatorial stochastic processes,
 for example, the Chinese restaurant process \cite{Pitman2005} and the
Beta-splitting model for fragmentation trees \cite{Aldous1996cladograms,McCullaghPitmanWinkel2008}. We need not prove
consistency explicitly for this model; it follows from our more general
construction of rewiring processes, all of which are consistent
(Theorem~\ref{thm:consistent rewiring chain}). We discuss the above
family further in Section~\ref{section:ER}.

%s3.2 #&#
\subsection{A more general construction}

Throughout the paper, we construct exchangeable and consistent Markov
processes using a special {\em rewiring measure} (Section~\ref{section:exchangeable rewiring maps}). In continuous-time, Markov
processes can admit infinitely many jumps in arbitrarily small time
intervals; however, by the consistency assumption, any edge can change
only finitely often in bounded intervals. In this case, we can choose a
$\sigma$-finite rewiring measure to direct the process.

%s4 #&#
\section{Preliminaries and the rewiring maps}\label{section:rewiring maps}
For $n=1,2,\ldots$\,, an (undirected) {\em graph} $G$ with vertex set
$[n]$ can be represented by its symmetric {\em adjacency matrix}
$(G_{ij})_{1\leq i,j\leq n}$ for which $G_{ij}=1$ if $G$ has an edge
between $i$ and $j$, and $G_{ij}=0$ otherwise.
By convention, we always assume $G_{ii}=0$ for all $i=1,\ldots,n$. We
write $ \mathcal{G}_n  $ to denote the finite
collection of all graphs with vertex set $[n]$.

On $\mathcal{G}_n  $, we define the following
operation of {\em rewiring}.
Let $w:=(w_{ij})_{1\leq i,j\leq n}$ be an $n\times n$ symmetric matrix
with entries in $\{0,1\}\times\{0,1\}$ and all diagonal entries $(0,0)$.
For convenience, we write each entry of $w$ as a pair
$w_{ij}:=(w_{ij}^0,w_{ij}^1)$, $1\leq i,j\leq n$.
We\vadjust{\goodbreak} define a map $w\dvtx  \mathcal{G}_n   \rightarrow
\mathcal{G}_n   $ by $G\mapsto G':=w(G)$, where
%
%e4.1 #&#
\begin{equation}
\label{eq:rewire} G'_{ij}:=\cases{ w_{ij}^0,&\quad
$G_{ij}=0$,
\cr
w_{ij}^1,&\quad $G_{ij}=1$, }\quad\quad
1\leq i,j\leq n.
\end{equation}
More compactly, we may write $w(G):=(w_{ij}^{G_{ij}})_{1\leq i,j\leq
n}$. We call $w$ a {\em rewiring map} and $w(G)$ the {\em rewiring} of
$G$ by $w$. We write $ \mathcal{W}_n   $ to denote
the collection of all rewiring maps $\mathcal{G}_n
\rightarrow \mathcal{G}_n  $, which are in
one-to-one correspondence with $n\times n$ symmetric matrices with
entries in $\{0,1\}\times\{0,1\}$ and all diagonal entries $(0,0)$.

The following display illustrates the rewiring operation in \eqref{eq:rewire}.
Given $G\in \mathcal{G}_n   $ and $w\in \mathcal{W}_n   $, we obtain $w(G)$ by choosing the
appropriate element of each entry of $w$: if $G_{ij}=0$, we choose the
{\em left} coordinate of $w_{ij}$; if $G_{ij}=1$, we choose the {\em
right} coordinate of $w_{ij}$.
For example,
\[
\begin{array} {c@{\quad\quad}ccc} G &  w & & w(G)
\\
\pmatrix{ 0 & 1 & 1 & 0 & 1
\cr
1 & 0 & 0 & 0 & 1
\cr
1 & 0 & 0 & 1 & 0
\cr
0 & 0 &
1 & 0 & 0
\cr
1 & 1 & 0 & 0 & 0 }&  \pmatrix{ ({\bf\underline0},0) & (1,{\bf
\underline0}) & (0,{\bf\underline 1}) & ({\bf\underline0},0) & (0,{\bf
\underline1})
\cr
(1,{\bf\underline0}) & ({\bf\underline0},0) & ({\bf\underline
1},0) & ({\bf\underline1},1) & (1,{\bf\underline0})
\cr
(0,{\bf\underline1}) & ({
\bf\underline1},0) & ({\bf\underline 0},0) & (0,{\bf\underline1}) & ({\bf
\underline0},0)
\cr
({\bf\underline0},0) & ({\bf\underline1},1) & (0,{\bf\underline
1}) & ({\bf\underline0},0) & ({\bf\underline1},0)
\cr
(0,{\bf\underline1}) & (1,{
\bf\underline0}) & ({\bf\underline 0},0) & ({\bf\underline1},0) & ({\bf
\underline0},0) } &\mapsto& \pmatrix{ 0 & 0 & 1 & 0 & 1
\cr
0 & 0 & 1 & 1 & 0
\cr
1
& 1 & 0 & 1 & 0
\cr
0 & 1 & 1 & 0 & 1
\cr
1 & 0 & 0 & 1 & 0 }. \end{array}
\]

A unique symmetric $n\times n$ matrix determines each element in
$ \mathcal{G}_n  $ and $ \mathcal
{W}_n  $, and so there is a natural restriction operation on
both spaces by taking the leading $m\times m$ submatrix, for any $m\leq
n$. In particular, we write
%
%e4.2 #&#
\begin{eqnarray}
\label{eq:restriction-graph} \Rmn G&:=&G|_{[m]}:=
(G_{ij})_{1\leq i,j\leq m},\quad\quad G
\in \mathcal{G}_n , \quad\mbox{and}\nonumber
\\[-8pt]\\[-8pt]
w|_{[m]}&:=&(w_{ij})_{1\leq i,j\leq m},\quad\quad w\in \mathcal
{W}_n ,
\nonumber
\end{eqnarray}
to denote the restrictions of $G\in \mathcal{G}_n
$ and $w\in \mathcal{W}_n  $ to $\graphsm$ and
$\wirem$, respectively.
These restriction operations lead to the notions of {\em infinite
graphs} and {\em infinite rewiring maps} as infinite symmetric arrays
with entries in the appropriate space, either $\{0,1\}$ or $\{0,1\}
\times\{0,1\}$. We write $\graphsN$ to denote the space of infinite
graphs, identified by a $\{0,1\}$-valued {\em adjacency array}, and
$\wireN$ to denote the space of infinite rewiring maps, identified by
a symmetric $\{0,1\}\times\{0,1\}$-valued array with $(0,0)$ on the diagonal.

Any $w\in\wireN$ acts on $\graphsN$ just as in \eqref{eq:rewire}
and, for any $G\in\graphsN$, the rewiring operation satisfies
\[
w(G)|_{[n]}=w|_{[n]}(G|_{[n]}) \quad\quad\mbox{for all }n\in
\mathbb{N}.
\]

The spaces $\graphsN$ and $\wireN$ are uncountable but can be
equipped with the discrete $\sigma$-algebras $\sigma \langle
\bigcup_{n\in\mathbb{N}} \mathcal{G}_n
\rangle$ and $\sigma \langle\bigcup_{n\in\mathbb{N}} \mathcal{W}_n   \rangle$, respectively, so that the
restriction maps $\cdot|_{[n]}$ are measurable for every $n\in\mathbb
{N}$. Moreover, both $\graphsN$ and $\wireN$ come equipped with a
product-discrete topology induced, for example, by the ultrametric
%
%e4.3 #&#
\begin{equation}
\label{eq:metric} d\bigl(w,w'\bigr):=1/\max\bigl\{n\in
\mathbb{N}\dvt w|_{[n]}=w'|_{[n]}\bigr\},\quad\quad
w,w'\in \wireN.
\end{equation}
The metric on $\graphsN$ is analogous. Both $\graphsN$ and $\wireN$
are compact, complete, and separable metric spaces.
Much of our development hinges on the following proposition, whose
proof is straightforward.\vspace*{-1pt}

%pr4.1 #&#
\begin{prop}\label{prop:Lipschitz}
Rewiring maps are associative under composition and Lipschitz
continuous in the metric \eqref{eq:metric}, with Lipschitz constant 1.\vspace*{-1pt}
\end{prop}

%s4.1 #&#
\subsection{Weakly exchangeable arrays}\label{section:arrays}
Let $\symmetricN$ denote the collection of {\em finite} permutations
of $\mathbb{N}$, that is, permutations $\sigma\dvtx \mathbb{N}\rightarrow
\mathbb{N}$ for which $\#\{i\in\mathbb{N}\dvt \sigma(i)\neq i\}<\infty
$. We call any random array $X:=(X_{ij})_{i,j\geq1}$ {\em weakly
exchangeable} if\vspace*{-1pt}
\[
X\mbox{ is almost surely symmetric, that is,
}X_{ij}=X_{ji}\mbox{ for all }i,j\mbox{ with probability
one},
\]
and
\[
X\equalinlaw X^{\sigma}:=(X_{\sigma(i)\sigma(j)})_{i,j\geq1}
\quad\quad\mbox{for all finite permutations }\sigma\dvtx \mathbb{N}\rightarrow \mathbb{N},
\]
where $\equalinlaw$ denotes {\em equality in law}.
Aldous defines weak exchangeability using only the latter condition;
see \cite{AldousExchangeability}, Chapter~14, page~132.
We impose symmetry for convenience~-- in this paper, all graphs and
rewiring maps are symmetric arrays.

From the discussion in Section~\ref{section:exchangeability}, we are
interested in models for random graphs $\Gamma$ that are {\em
exchangeable}, meaning the adjacency matrix $(\Gamma_{ij})_{i,j\geq
1}$ is a weakly exchangeable $\{0,1\}$-valued array. Likewise, we call
a random rewiring map $W$ {\em exchangeable} if its associated $\{0,1\}
\times\{0,1\}$-valued array $(W_{ij})_{i,j\geq1}$ is weakly exchangeable.

de Finetti's theorem represents any infinitely exchangeable sequence
$Z:=(Z_i)_{i\geq1}$ in a Polish space $\mathcal{S}$ with a
(non-unique) measurable function $g\dvtx [0,1]^2\rightarrow\mathcal{S}$
such that $Z\equalinlaw Z^*$, where
%
%e4.4 #&#
\begin{equation}
\label{eq:de Finetti rep}Z_i^*:=g(\alpha,\eta _i),\quad\quad i\geq1,
\end{equation}
for $\{\alpha;(\eta_i)_{i\geq1}\}$ independent, identically
distributed (i.i.d.)~Uniform random variables on $[0,1]$.
The Aldous--Hoover theorem \cite{AldousExchArrays,AldousExchangeability} extends de Finetti's
representation \eqref{eq:de Finetti rep} to weakly exchangeable
$\mathcal{S}$-valued arrays: to any such array $X$, there exists a
(non-unique) measurable function $f\dvtx [0,1]^4\rightarrow\mathcal{S}$
satisfying $f(\bullet,b,c,\bullet)=f(\bullet,c,b,\bullet)$ such
that $X\equalinlaw X^*$, where
%
%e4.5 #&#
\begin{equation}
\label{eq:A-H rep} X^*_{ij}:=f(\alpha,\eta_i,
\eta_j,\lambda_{\{i,j\}}),\quad\quad i,j\geq1,
\end{equation}
for $\{\alpha;(\eta_{i})_{i\geq1};(\lambda_{\{i,j\}})_{i>j\geq1}\}
$ i.i.d. Uniform random variables on $[0,1]$.

The function $f$ has a statistical interpretation that reflects the
structure of the random array.
In particular, $f$ decomposes the law of $X^*_{ij}$ into individual
$\lambda_{\{i,j\}}$, row $\eta_i$, column $\eta_j$, and overall
$\alpha$ effects.
The overall effect plays the role of the mixing measure in the de
Finetti interpretation.
If $g$ in \eqref{eq:de Finetti rep} is constant with respect to its
first argument,  that is,  $g(a,\cdot)=g(a',\cdot)$ for all $a,a'\in
[0,1]$, then $Z^*$ constructed in \eqref{eq:de Finetti rep} is an
i.i.d. sequence. Letting $g$ vary with its first argument produces a
mixture of i.i.d. sequences.
A fundamental interpretation of de Finetti's theorem is:
\[
\mbox{every infinitely exchangeable sequence is a mixture of i.i.d. sequences.}
\]
Similarly, if $f$ in \eqref{eq:A-H rep} satisfies $f(a,\cdot,\cdot
,\cdot)=f(a',\cdot,\cdot,\cdot)$ for all $a,a'\in[0,1]$, then
$X^*$ is {\em dissociated}, that is
%
%e4.6 #&#
\begin{equation}
\label{eq:dissociated} X^*|_{[n]}\mbox{ is independent of }X^*|_{\{n+1,n+2,\ldots\}}
\mbox{ for all }n\in\mathbb{N}.
\end{equation}
The Aldous--Hoover representation \eqref{eq:A-H rep} spurs the sequel
to de Finetti's interpretation:
\[
\mbox{every weakly exchangeable array is a mixture of dissociated arrays.}
\]
See Aldous \cite{AldousExchangeability}, Chapter~14, for more details.
We revisit the theory of weakly exchangeable arrays in Section~\ref{section:exchangeable rewiring maps}.

%s5 #&#
\section{Discrete-time rewiring Markov chains}\label{section:discrete}
Throughout the paper, we use the rewiring maps to construct Markov
chains on $\graphsN$. From any probability distribution $\omega_n$ on
$ \mathcal{W}_n  $, we generate $W_1,W_2,\ldots$
i.i.d. from $\omega_n$ and a random graph $\Gamma_0\in
\mathcal{G}_n  $ (independently of $W_1,W_2,\ldots$). We
then define a Markov chain $\{\Gamma_m\}_{m=0,1,2,\ldots}$ on
$ \mathcal{G}_n  $ by
%
%e5.1 #&#
\begin{equation}
\label{eq:rewiring chain} \Gamma_{m}:=W_{m}(\Gamma_{m-1})=(W_m
\circ\cdots\circ W_1) (\Gamma _0),\quad\quad m\geq1.
\end{equation}
We call $\omega_n$ {\em exchangeable} if $W\sim\omega_n$ is an
exchangeable rewiring map, that is, $W\equalinlaw W^{\sigma}$ for all
permutations $\sigma\dvtx [n]\rightarrow[n]$.

%pr5.1 #&#
\begin{prop}\label{prop:exchangeable rewiring chain}
Let $\omega_n$ be an exchangeable probability measure on $ \mathcal{W}_n  $ and let $\boldsymbol{\Gamma}:=\{\Gamma_m\}
_{m=0,1,2,\ldots}$ be as constructed in \eqref{eq:rewiring chain}
from an exchangeable initial state $\Gamma_0$ and $W_1,W_2,\ldots$
i.i.d. from $\omega_n$. Then $\boldsymbol{\Gamma}$ is an exchangeable
Markov chain on $ \mathcal{G}_n  $ with transition
probability
%
%e5.2 #&#
\begin{equation}
\label{eq:rewiring tps} P_{\omega_n}\bigl(G,G'\bigr):=
\omega_n\bigl(\bigl\{W\in \mathcal {W}_n
\dvt W(G)=G'\bigr\}\bigr),\quad\quad G,G'\in \mathcal
{G}_n .
\end{equation}
\end{prop}

\begin{pf}
The Markov property is immediate by mutual independence of $\Gamma
_0,W_1,W_2,\ldots$\,. The formula for the transition probabilities
\eqref{eq:rewiring tps} follows by description \eqref{eq:rewiring
chain} of $\boldsymbol{\Gamma}$.

We need only show that $\boldsymbol{\Gamma}$ is exchangeable. By
assumption, $\Gamma_0$ is an exchangeable random graph on $n$
vertices, and so its distribution is invariant under arbitrary
permutation of $[n]$. Moreover, the law of $W\sim\omega_n$ satisfies
$W\equalinlaw W^{\sigma}$ and, for any fixed $w\in \mathcal
{W}_n  $ and $\sigma\in \mathscr{S}_{n}
$, $G':=w(G)$ satisfies
\[
G'^{\sigma}_{ij}:=G'_{\sigma(i)\sigma(j)}=w_{\sigma(i)\sigma
(j)}^{G_{\sigma(i)\sigma(j)}},
\]
the $ij$-entry of $w^{\sigma}(G^{\sigma})$.
Therefore, $W^{\sigma}(G^{\sigma})=W(G)^{\sigma}$ and, for any
exchangeable graph $\Gamma$ and exchangeable rewiring map $W$, we have
\[
W(\Gamma)^{\sigma}=W^{\sigma}\bigl(\Gamma^{\sigma}\bigr)
\equalinlaw W(\Gamma ) \quad\quad\mbox{for all }\sigma\in \mathscr{S}_{n} .
\]
Hence, the transition law of $\boldsymbol{\Gamma}$ is equivariant with
respect to relabeling. Since the initial state $\Gamma_0$ is
exchangeable, so is the Markov chain.
\end{pf}
%
%de5.1 #&#
\begin{defn}We call $\{\Gamma_m\}_{m=0,1,2,\ldots}$ an {\em$\omega
_n$-rewiring Markov chain}.
\end{defn}

From the discussion in Section~\ref{section:rewiring maps}, we can
define an exchangeable measure $\omega^{(n)}$ on $ \mathcal
{W}_n  $ as the restriction to $ \mathcal
{W}_n  $ of an exchangeable probability measure $\omega$ on
$\wireN$, where
%
%e5.3 #&#
\begin{equation}
\label{eq:omega-n omega}\omega^{(n)}(W):=\omega\bigl(\bigl\{ W^*\in
\wireN\dvt W^*|_{[n]}=W\bigr\}\bigr),\quad\quad W\in \mathcal {W}_n .
\end{equation}
Denote by $P_{\omega}^{(n)}$ the transition probability measure of an
$\omega^{(n)}$-rewiring Markov chain on $ \mathcal
{G}_n  $, as defined in \eqref{eq:rewiring tps}.

%th5.1 #&#
\begin{theorem}\label{thm:consistent rewiring chain}
For any exchangeable probability measure $\omega$ on $\wireN$, $\{
P_{\omega}^{(n)}\}_{n\in\mathbb{N}}$ is a consistent family of
exchangeable transition probabilities in the sense that
%
%e5.4 #&#
\begin{equation}
\label{eq:consistent tps} P_{\omega}^{(m)}\bigl(G,G'
\bigr)=P_{\omega}^{(n)}\bigl(G^*,\mathbf {R}^{-1}_{m,n}
\bigl(G'\bigr)\bigr),\quad\quad G,G'\in\graphsm, \quad\quad\mbox{for all }m
\leq n,
\end{equation}
for every $G^*\in\mathbf{R}_{m,n}^{-1}(G):=\{G''\in \mathcal
{G}_n  \dvt G''|_{[m]}=G\}$, where $\Rmn$ is defined in \eqref
{eq:restriction-graph}.
\end{theorem}

\begin{pf}
Proposition~\ref{prop:exchangeable rewiring chain} implies
exchangeability of $\{P_{\omega}^{(n)}\}_{n\in\mathbb{N}}$. It
remains to show that $\{P_{\omega}^{(n)}\}_{n\in\mathbb{N}}$
satisfies \eqref{eq:consistent tps}.
By \eqref{eq:rewiring tps},
\[
P_{\omega}^{(n)}\bigl(G,G'\bigr):=
\omega^{(n)}\bigl(\bigl\{W\in \mathcal {W}_n
\dvt W(G)=G'\bigr\}\bigr),\quad\quad G,G'\in \mathcal
{G}_n .
\]
Now, for any $m\leq n$, fix $G,G'\in\graphsm$ and $G^*\in\mathbf
{R}^{-1}_{m,n}(G)$. Then \eqref{eq:consistent tps} requires
\[
\omega^{(n)}\bigl(\bigl\{W\in \mathcal{W}_n \dvt W\bigl(G^*
\bigr)\in \mathbf{R}^{-1}_{m,n}\bigl(G'\bigr)\bigr
\}\bigr)=\omega^{(m)}\bigl(\bigl\{W\in\wirem\dvt W(G)=G'\bigr\}
\bigr),
\]
which follows by definition \eqref{eq:omega-n omega} of $\omega
^{(n)}$. To see this, note that
\begin{eqnarray*}
\omega^{(n)}\bigl(\bigl\{W\in \mathcal{W}_n
\dvt W|_{[m]}(G)=G'\bigr\}\bigr)&=&\omega\bigl(\bigl\{W\in
\wireN\dvt (W|_{[n]})|_{[m]}(G)=G'\bigr\}\bigr)
\\
&=&\omega\bigl(\bigl\{W\in\wireN\dvt W|_{[m]}(G)=G'\bigr\}
\bigr)
\\
&=&\omega^{(m)}\bigl(\bigl\{W\in\wirem\dvt W(G)=G'\bigr\}
\bigr).
\end{eqnarray*}
This completes the proof.
\end{pf}

%re5.1 #&#
\begin{rmk}
The consistency condition \eqref{eq:consistent tps} for Markov chains
is exactly the necessary and sufficient condition for a function of a
Markov chain to be a Markov chain, as proven in \cite{BurkeRosenblatt1958}. Before describing the measure $\omega$ from
Theorem~\ref{thm:consistent rewiring chain} in further detail, we
first show some concrete examples of rewiring chains.
\end{rmk}

%s5.1 #&#
\subsection{The Erd\H{o}s--R\'enyi rewiring chain}\label{section:ER}

For any $0\leq p\leq1$, let $\varepsilon_p$ denote the Erd\H{o}s--R\'
enyi measure on $\graphsN$, which we define by its finite-dimensional
restrictions $\varepsilon^{(n)}_p$ to $ \mathcal
{G}_n  $ for each $n\in\mathbb{N}$,
%
%e5.5 #&#
\begin{equation}
\label{eq:ER fidi} \varepsilon^{(n)}_p(G):=\prod
_{1\leq i<j\leq
n}p^{G_{ij}}(1-p)^{1-G_{ij}},\quad\quad G\in \mathcal
{G}_n .
\end{equation}
Given any pair $(p_0,p_1)\in[0,1]\times[0,1]$, the {\em
$(p_0,p_1)$-Erd\H{o}s--R\'enyi chain} has finite-dimensional transition
probabilities
%
%e5.6 #&#
\begin{equation}
\label{eq:ER fidi tps} P^{(n)}_{p_0,p_1}\bigl(G,G'\bigr):=
\prod_{1\leq i<j\leq
n}p_{G_{ij}}^{G'_{ij}}(1-p_{G_{ij}})^{1-G'_{ij}},
\quad\quad G,G'\in \mathcal{G}_n .
\end{equation}

%pr5.2 #&#
\begin{prop}
For $0<p_0,p_1<1$, the $(p_0,p_1)$-Erd\H{o}s--R\'enyi rewiring chain
has unique stationary distribution $\varepsilon_q$, with $q:=p_0/(1-p_1+p_0)$.
\end{prop}
\begin{pf}
By assumption, both $p_0$ and $p_1$ are strictly between $0$ and $1$
and, thus, \eqref{eq:ER fidi} assigns positive probability to every
transition in $\mathcal{G}_n  $, for every $n\in
\mathbb{N}$. Therefore, each finite-dimensional chain is aperiodic and
irreducible, and each possesses a unique stationary distribution
$\theta^{(n)}$.
By consistency of the transition probabilities $\{P_{p_0,p_1}^{(n)}\}
_{n\in\mathbb{N}}$ (Theorem~\ref{thm:consistent rewiring chain}),
the finite-dimensional stationary measures $\{\theta^{(n)}\}_{n\in
\mathbb{N}}$ must be exchangeable and consistent and, therefore, they
determine a unique measure $\theta$ on $\graphsN$, which is
stationary for $P_{p_0,p_1}$. Furthermore, by conditional independence
of the edges of $G'$, given $G$, the stationary law must be Erd\H
{o}s--R\'enyi with some parameter $q\in(0,1)$.

In an $\varepsilon_q^{(n)}$-random graph, all edges are present or not
independently with probability $q$. Therefore, it suffices to look at
the probability of the edge between vertices labeled 1 and 2. In this
case, we need to choose $q$ so that
\[
qp_1+(1-q)p_0=q,
\]
which implies $q=p_0/(1-p_1+p_0)$.
\end{pf}

%re5.2 #&#
\begin{rmk}
Some elementary special cases of the $(p_0,p_1)$-Erd\H{o}s--R\'enyi
rewiring chain are worth noting. First, for $(p_0,p_1)$ either $(0,0)$
or $(1,1)$, this chain is degenerate at either the empty graph $\mathbf
{0}$ or the complete graph $\mathbf{1}$ and has unique stationary
measure $\varepsilon_0$ or $\varepsilon_1$, respectively. On the
other hand, when $(p_0,p_1)=(0,1)$, the chain is degenerate at its
initial state and so its initial distribution is stationary. However,
if $(p_0,p_1)=(1,0)$, then the chain alternates between its initial
state $G$ and its complement $\bar{G}:=(\bar{G}_{ij})_{i,j\geq1}$,
where $\bar{G}_{ij}:=1-G_{ij}$ for all $i,j\geq1$; in this case, the
chain is periodic and does not have a unique stationary distribution.
We also note that when $(p_0,p_1)=(p,p)$ for some $p\in(0,1)$, the
chain is simply an i.i.d. sequence of $\varepsilon_p$-random graphs
with stationary distribution $\varepsilon_q$, where $q=p/(1-p+p)=p$,
as it must.
\end{rmk}

For $\alpha,\beta>0$, we define the {\em mixed Erd\H{o}s--R\'enyi
rewiring chain} through $\varepsilon_{\alpha,\beta}^{(n)}$, the
mixture of $\varepsilon_p^{(n)}$-laws with respect to the Beta law
with parameter $(\alpha,\beta)$.
Writing
\[
\mathscr{B}_{\alpha,\beta}(\mathrm{d}p)=\frac{\Gamma(\alpha+\beta
)}{\Gamma(\alpha)\Gamma(\beta)}p^{\alpha-1}(1-p)^{\beta-1}\,\mathrm{d}p,\vadjust{\goodbreak}
\]
we derive\vspace*{-1pt}
\begin{eqnarray*}
\varepsilon_{\alpha,\beta}^{(n)}(G)&:=&\int_{[0,1]}
\varepsilon _p^{(n)}(G)\frac{\Gamma(\alpha+\beta)}{\Gamma(\alpha)\Gamma
(\beta)}p^{\alpha-1}(1-p)^{\beta-1}\,\mathrm{d}p
\\[-1pt]
&=&\frac{\Gamma(\alpha+\beta)}{\Gamma(\alpha)\Gamma(\beta
)}\frac{\Gamma(\alpha+n_1)\Gamma(\beta+n_0)}{\Gamma(\alpha+\beta
+n)}\int_{[0,1]}
\mathscr{B}_{\alpha+n_1,\beta+n_0}(\mathrm{d}p)
\\[-1pt]
&=&\frac{\alpha^{\uparrow n_1}\beta^{\uparrow n_0}}{(\alpha+\beta
)^{\uparrow n}},
\end{eqnarray*}
where $n_r:=\sum_{1\leq i<j\leq n}\{G_{ij}=r\}$, $r=0,1$, and
$\alpha^{\uparrow n}=\alpha(\alpha+1)\cdots(\alpha+n-1)$.
For $\alpha_0,\beta_0,\alpha_1,\beta_1>0$, we define {\em mixed
Erd\H{o}s--R\'enyi transition probabilities} by\vspace*{-1pt}
%
%e5.7 #&#
\begin{equation}
\label{eq:mixture ER} P^{(n)}_{(\alpha_0,\beta_0),(\alpha_1,\beta_1)}\bigl(G,G'\bigr):=
\frac
{\alpha_0^{\uparrow n_{01}}\beta_0^{\uparrow n_{00}}\alpha
_1^{\uparrow n_{11}}\beta_1^{\uparrow n_{10}}}{(\alpha_0+\beta
_0)^{\uparrow(n_{00}+n_{01})}(\alpha_1+\beta_1)^{\uparrow
(n_{10}+n_{11})}},\quad\quad G,G'\in\mathcal{G}_n .
\end{equation}

An interesting special case takes $(\alpha_0,\beta_0)=(\beta,\alpha
)$ and $(\alpha_1,\beta_1)=(\alpha',\beta)$ for $\alpha,\alpha
',\beta>0$. In this case, \eqref{eq:mixture ER} becomes\vspace*{-1pt}
\[
P^{(n)}_{(\beta,\alpha),(\alpha',\beta)}\bigl(G,G'\bigr)=\frac{\alpha
^{\uparrow n_{00}}\beta^{\uparrow n_{01}}\alpha'^{\uparrow
n_{11}}\beta^{\uparrow n_{10}}}{(\alpha+\beta)^{\uparrow n_0}(\alpha
'+\beta)^{\uparrow n_1}},\quad\quad
G,G'\in\mathcal{G}_n .
\]

%pr5.3 #&#
\begin{prop}
$P^{(n)}_{(\beta,\alpha),(\alpha',\beta)}$ is reversible with
respect to $\varepsilon_{\alpha+\beta,\alpha'+\beta}^{(n)}$.\vspace*{-2pt}
\end{prop}
\begin{pf}
For fixed $G,G'\in\mathcal{G}_n  $, we write
$n_{rs}:=\sum_{i<j}\mathbf{1}\{G_{ij}=r\mbox{ and }G'_{ij}=s\}$ and
$n'_{rs}:=\sum_{i<j}\mathbf{1}\{G'_{ij}=r\mbox{ and }G_{ij}=s\}$.
Note that $n'_{rs}=n_{sr}$. Therefore, we have\vspace*{-1pt}
\begin{eqnarray*}
\varepsilon^{(n)}_{\alpha+\beta,\alpha'+\beta}(G)P_{(\beta,\alpha
),(\alpha',\beta)}^{(n)}
\bigl(G,G'\bigr)&=&\frac{\alpha^{\uparrow
n_{00}}\beta^{\uparrow n_{01}}\alpha'^{\uparrow n_{11}}\beta
^{\uparrow n_{10}}}{(\alpha+2\beta+\alpha')^{\uparrow n}}
\\[-1pt]
&=&\frac{\alpha^{\uparrow n'_{00}}\beta^{\uparrow n'_{10}}\alpha
'^{\uparrow n'_{11}}\beta^{\uparrow n'_{01}}}{(\alpha+2\beta+\alpha
')^{\uparrow n}}
\\[-1pt]
&=&\varepsilon^{(n)}_{\alpha+\beta,\alpha'+\beta}\bigl(G'
\bigr)P_{(\beta
,\alpha),(\alpha',\beta)}^{(n)}\bigl(G',G\bigr),
\end{eqnarray*}
establishing detailed balance and, thus, reversibility.\vspace*{-1pt}
\end{pf}

A mixed Erd\H{o}s--R\'enyi Markov chain is directed by\vspace*{-1pt}
\[
\omega(\mathrm{d}W):=\int_{[0,1]\times[0,1]}\omega_{p_0,p_1}(\mathrm{d}W) (\mathscr
{B}_{\alpha_0,\beta_0}\otimes\mathscr{B}_{\alpha_1,\beta
_1}) (\mathrm{d}p_0,\mathrm{d}p_1),\quad\quad
W\in\wireN,
\]
where $\omega_{p_0,p_1}$ is determined by its finite-dimensional distributions\vspace*{-1pt}
\[
\omega_{p_0,p_1}^{(n)}(W):=\prod_{1\leq i<j\leq
n}p_0^{W^0_{ij}}(1-p_0)^{1-W_{ij}^0}p_1^{W_{ij}^1}(1-p_1)^{1-W_{ij}^1},\quad\quad
W\in \mathcal{W}_n ,
\]
for $0<p_0,p_1<1$, for every $n\in\mathbb{N}$.\vadjust{\goodbreak}

In the next section, we see that a representation of the directing
measure $\omega$ as a mixture of simpler measures holds more
generally. Notice that $W\sim\omega_{p_0,p_1}$ is dissociated for all
fixed $(p_0,p_1)\in(0,1)\times(0,1)$. By the Aldous--Hoover theorem,
we can express any exchangeable measure on $\wireN$ as a mixture of
dissociated measures.

%s6 #&#
\section{Exchangeable rewiring maps and their rewiring limits}\label
{section:exchangeable rewiring maps}
To more precisely describe the mixing measure $\omega$, we extend the
theory of {\em graph limits} to its natural analog for rewiring maps.
We first review the related theory of graph limits, as surveyed by Lov\'
asz \cite{LovaszBook}.

%s6.1 #&#
\subsection{Graph limits}\label{section:graph limits}
A graph limit is a statistic that encodes a lot of structural
information about an infinite graph.
In essence, the graph limit of an exchangeable random graph contains
all relevant information about its distribution.

For any injection $\psi\dvtx [m]\rightarrow[n]$, $m\leq n$, and $G\in
\mathcal{G}_n  $, we define $G^{\psi}:=(G_{\psi
(i)\psi(j)})_{1\leq i,j\leq m}$.
In words, $G^{\psi}$ is the subgraph $G$ induces on $[m]$ by the
vertices in the range of $\psi$. Given $G\in\mathcal{G}_n  $ and $F\in\graphsm$, we define
$\ind(F,G)$ to equal the number of injections $\psi\dvtx [m]\rightarrow
[n]$ such that $G^{\psi}=F$. Intuitively, $\ind(F,G)$ is the number
of ``copies'' of $F$ in $G$, which we normalize to obtain the {\em
density} of $F$ in $G$,
%
%e6.1 #&#
\begin{equation}
\label{eq:density} t(F,G):=\frac{\ind(F,G)}{n^{\downarrow m}},\quad\quad F\in\graphsm, G\in
\mathcal{G}_n ,
\end{equation}
where $n^{\downarrow m}:=n(n-1)\cdots(n-m+1)$ is the number of unique
injections $\psi\dvtx [m]\rightarrow[n]$. The {\em limiting density} of
$F$ in any infinite graph $G\in\graphsN$ is
%
%e6.2 #&#
\begin{equation}
\label{eq:limiting density} t(F,G):=\lim_{n\rightarrow\infty}t(F,G|_{[n]}),\quad\quad F\in
\graphsm , \quad\quad\mbox{if it exists}.
\end{equation}

The collection $\mathcal{G}^*:=\bigcup_{m\in\mathbb{N}}\graphsm$
is countable and so we can define the {\em graph limit} of $G\in
\graphsN$ by
%
%e6.3 #&#
\begin{equation}
\label{eq:graph limit} |G|:=\bigl(t(F,G)\bigr)_{F\in\mathcal{G}^*},
\end{equation}
provided $t(F,G)$ exists for all $F\in\mathcal{G}^*$. Any graph limit
is an element in $[0,1]^{\mathcal{G}^*}$, which is compact under the
metric
%
%e6.4 #&#
\begin{equation}
\label{eq:graph limit metric} \rho\bigl(x,x'\bigr):=\sum
_{n\in\mathbb{N}}2^{-n}\sum_{F\in\mathcal{G}_n  }\bigl|x_F-x'_F\bigr|,\quad\quad
x,x'\in[0,1]^{\mathcal{G}^*}.
\end{equation}
The space of graph limits is a compact subset of $[0,1]^{\mathcal
{G}^{*}}$, which we denote by $\graphlimits$. We implicitly equip
$[0,1]^{\mathcal{G}^*}$ with its Borel $\sigma$-field and
$\graphlimits$ with its trace $\sigma$-field.

Any $D\in\graphlimits$ is a sequence $(D_F)_{F\in\mathcal{G}^*}$,
where $D(F):=D_F$ denotes the coordinate of $D$ corresponding to $F\in
\mathcal{G}^*$. In this way, any $D\in\graphlimits$ determines a
probability measure $\gamma_D^{(n)}$ on $\mathcal{G}_n  $, for every $n\in\mathbb{N}$, by
%
%e6.5 #&#
\begin{equation}
\label{eq:graph limit induced measure} \gamma_D^{(n)}(G):=D(G),\quad\quad G\in
\mathcal{G}_n .
\end{equation}
Furthermore, the collection $(\gamma_D^{(n)})_{n\in\mathbb{N}}$ is
consistent and exchangeable on $\{\mathcal{G}_n  \}
_{n\in\mathbb{N}}$ and, by Kolmogorov's extension theorem, determines
a unique exchangeable measure $\gamma_D$ on $\graphsN$, for which
$\gamma_D$-almost every $G\in\graphsN$ has $|G|=D$.

Conversely, combining the Aldous--Hoover theorem for weakly exchangeable
arrays (\cite{AldousExchangeability}, Theorem~14.21) and Lov\'asz--Szegedy theorem of graph limits (\cite{LovaszSzegedy2006}, Theorem~2.7), any exchangeable random graph $\Gamma$ is governed by a mixture
of $\gamma_D$ measures. In particular, to any exchangeable random
graph $\Gamma$, there exists a unique probability measure $\Delta$ on
$\graphlimits$ such that $\Gamma\sim\gamma_{\Delta}$, where
%
%e6.6 #&#
\begin{equation}
\label{eq:Delta-mixture} \gamma_{\Delta}(\cdot):=\int_{\graphlimits}
\gamma_D(\cdot)\Delta(\mathrm{d}D).
\end{equation}

%s6.2 #&#
\subsection{Rewiring limits}\label{section:rewiring limits}
Since $\{0,1\}\times\{0,1\}$ is a finite space, the Aldous--Hoover
theorem applies to exchangeable rewiring maps. Following Section~\ref{section:graph limits}, we define the {\em density} of $V\in\wirem$
in $W\in \mathcal{W}_n  $ by
%
%e6.7 #&#
\begin{equation}
\label{eq:rewire density} t(V,W):=\frac{\ind(V,W)}{n^{\downarrow m}},
\end{equation}
where $\ind(V,W)$ equals the number of injections $\psi
\dvtx [m]\rightarrow[n]$ for which $W^{\psi}=V$. For an infinite rewiring
map $W\in\wireN$, we define
\[
t(V,W):=\lim_{n\rightarrow\infty}t(V,W|_{[n]}),\quad\quad W\in\mathcal
{W}_m,\quad\quad \mbox{if it exists}.
\]
As for graphs, the collection $\mathcal{W}^{*}:=\bigcup_{m\in\mathbb
{N}}\wirem$ is countable and so we can define the {\em rewiring limit}
of $W\in\wireN$ by
%
%e6.8 #&#
\begin{equation}
\label{eq:rewiring limit} |W|:=\bigl(t(V,W)\bigr)_{V\in\mathcal{W}^*},
\end{equation}
provided $t(V,W)$ exists for all $V\in\mathcal{W}^*$.

We write $\rewiringlimits\subset[0,1]^{\mathcal{W}^*}$ to denote the
compact space of rewiring limits and $\upsilon_V=\upsilon(V)$ to
denote the coordinate of $\upsilon\in\rewiringlimits$ corresponding
to $V\in\mathcal{W}^*$. We equip $\rewiringlimits$ with the metric
%
%e6.9 #&#
\begin{equation}
\label{eq:rewiring metric} \rho\bigl(\upsilon,\upsilon'\bigr):=\sum
_{n\in\mathbb{N}}2^{-n}\sum_{V\in
\mathcal{W}_n  }\bigl|
\upsilon_V-\upsilon'_V\bigr|,\quad\quad \upsilon,
\upsilon'\in\rewiringlimits.
\end{equation}

%le6.1 #&#
\begin{lemma}\label{lemma:rewiring structure}
Every $\upsilon\in\rewiringlimits$ satisfies
\begin{itemize}
\item$\upsilon(V)=\sum_{\{V^*\in\mathcal{W}_{n+1}:V^*|_{[n]}=V\}
}\upsilon(V^*)$ for every $V\in \mathcal{W}_n  $,
for all $n\in\mathbb{N}$, and
\item$\sum_{V\in \mathcal{W}_n  }\upsilon(V)=1$
for every $n\in\mathbb{N}$.
\end{itemize}
\end{lemma}

\begin{pf}
By definition of $\rewiringlimits$, we may assume that $\upsilon$ is
the rewiring limit $|W^*|$ of some $W^*\in\wireN$ so that $\upsilon
(V)=t(V,W^*)$, for every $V\in\mathcal{W}^*$.
From the definition of the rewiring limit \eqref{eq:rewiring limit},
\[
\sum_{W\in\wirem}\upsilon(W)=\sum
_{W\in\wirem}\lim_{n\rightarrow\infty}\frac{\ind(W,W^*|_{[n]})}{n^{\downarrow
m}}=\lim
_{n\rightarrow\infty}\sum_{W\in\wirem}
\frac{\ind
(W,W^*|_{[n]})}{n^{\downarrow m}}=1,
\]
where the interchange of sum and limit is justified by the Bounded
Convergence theorem because $0\leq\ind(W,W^*|_{[n]})/n^{\downarrow
m}\leq1$ for all $W\in\wirem$. Also, for every $m\leq n$ and $W\in
\wirem$, we have
\begin{eqnarray*}
\sum_{\{W'\in \mathcal{W}_n  :W'|_{[m]}=W\}
}\upsilon\bigl(W'\bigr)&=&
\sum_{\{W'\in \mathcal{W}_n
:W'|_{[m]}=W\}}\lim_{k\rightarrow\infty}
\frac{\ind
(W',W^*|_{[k]})}{k^{\downarrow n}}
\\
&=&\lim_{k\rightarrow\infty}\sum_{\{W'\in \mathcal
{W}_n  :W'|_{[m]}=W\}}
\frac{\ind
(W',W^*|_{[k]})}{k^{\downarrow n}}
\\
&=&\lim_{k\rightarrow\infty}\frac{\ind
(W,W^*|_{[k]})}{k^{\downarrow m}}
\\
&=&v(W).
\end{eqnarray*}
This follows by the definition of $\ind(\cdot,\cdot)$ and also
because, for any $\psi\dvtx [m]\rightarrow[k]$ there are $k^{\downarrow
n}/k^{\downarrow m}$ injections $\psi'\dvtx [n]\rightarrow[k]$ such that
$\psi'$ coincides with $\psi$ on $[m]$.
\end{pf}

%le6.2 #&#
\begin{lemma}\label{lemma:compact rewiring}
$(\rewiringlimits,\rho)$ is a compact metric space.
\end{lemma}

%th6.1 #&#
\begin{theorem}\label{prop:dissociated}
Let $W$ be a dissociated exchangeable rewiring map. Then, with
probability one, $|W|$ exists and is nonrandom.
\end{theorem}

We delay the proofs of Lemma~\ref{lemma:compact rewiring} and Theorem~\ref{prop:dissociated} until Section~\ref{section:proof}.

%co6.1 #&#
\begin{cor}\label{cor:existence rewiring limit}
Let $W\in\wireN$ be an exchangeable random rewiring map. Then $|W|$
exists almost surely.
\end{cor}

\begin{pf}
By Theorem~\ref{prop:dissociated}, every dissociated rewiring map
possesses a nonrandom rewiring limit almost surely.
By the Aldous--Hoover theorem, $W$ is a mixture of dissociated rewiring
maps and the conclusion follows.
\end{pf}

By Lemma~\ref{lemma:rewiring structure}, any $\upsilon\in
\rewiringlimits$ determines a probability measure $\Omega_{\upsilon
}$ on $\wireN$ in a straightforward way: for each $n\in\mathbb{N}$,
we define $\Omega_{\upsilon}^{(n)}$ as the probability distribution
on $ \mathcal{W}_n  $ with
%
%e6.10 #&#
\begin{equation}
\label{eq:measure on wiren} \Omega_{\upsilon}^{(n)}(w):=\upsilon(w),
\quad\quad w\in
\mathcal {W}_n .
\end{equation}

%pr6.1 #&#
\begin{prop}\label{prop:rewiring measure}
For any $\upsilon\in\rewiringlimits$, $\{\Omega_{\upsilon}^{(n)}\}
_{n\in\mathbb{N}}$ is a collection of exchangeable and consistent
probability distributions on $\{ \mathcal{W}_n  \}
_{n\in\mathbb{N}}$. In particular, $\{\Omega_{\upsilon}^{(n)}\}
_{n\in\mathbb{N}}$ determines a unique exchangeable probability
measure $\Omega_{\upsilon}$ on $\wireN$ for which $\Omega_{\upsilon
}$-almost every $w\in\wireN$ has $|w|=\upsilon$.
\end{prop}

\begin{pf}
By Lemma~\ref{lemma:rewiring structure}, the collection $\{\Omega
_{\upsilon}^{(n)}\}_{n\in\mathbb{N}}$ in \eqref{eq:measure on
wiren} is a consistent family of probability distributions on $\{
\mathcal{W}_n  \}_{n\in\mathbb{N}}$.
Exchangeability follows because $\ind(w,W^*|_{[n]})$ is invariant
under relabeling of $w$, that is, $\ind(w,W^*)=\ind(w^{\sigma
},W^*|_{[n]}^{\sigma'})$ for all permutations $\sigma\in\mathscr
{S}_m$ and $\sigma'\in\mathscr{S}_{n} $. By
Kolmogorov's extension theorem, $\{\Omega_{\upsilon}^{(n)}\}_{n\in
\mathbb{N}}$ determines a unique measure $\Omega_{\upsilon}$ on the
limit space $\wireN$. Finally, $W\sim\Omega_{\upsilon}$ is
dissociated and so, by Theorem~\ref{prop:dissociated}, $|W|=\upsilon$
almost surely.
\end{pf}

We call $\Omega_{\upsilon}$ in Proposition~\ref{prop:rewiring
measure} a {\em rewiring measure} directed by $\upsilon$. For any
measure $\Upsilon$ on $\rewiringlimits$, we define the $\Upsilon
$-mixture of rewiring measures by
%
%e6.11 #&#
\begin{equation}
\label{eq:rewiring mixture} \Omega_{\Upsilon}(\cdot):=\int_{\rewiringlimits}
\Omega_{\upsilon
}(\cdot)\Upsilon(\mathrm{d}\upsilon).
\end{equation}

%co6.2 #&#
\begin{cor}
To any exchangeable rewiring map $W$, there exists a unique probability
measure $\Upsilon$ on $\rewiringlimits$ such that $W\sim\Omega
_{\Upsilon}$.
\end{cor}
\begin{pf}
This follows by the Aldous--Hoover theorem and Proposition~\ref
{prop:rewiring measure}.
\end{pf}

From Theorem~\ref{prop:dissociated} and Proposition~\ref
{prop:rewiring measure}, any probability measure $\Upsilon$ on
$\rewiringlimits$ corresponds to an $\Omega_{\Upsilon}$-rewiring
chain as in Theorem~\ref{thm:consistent rewiring chain}.

%s7 #&#
\section{Continuous-time rewiring processes}\label{section:continuous}
We now refine our discussion to rewiring chains in continuous-time, for
which infinitely many transitions can ``bunch up'' in arbitrarily small
intervals, but individual edges jump only finitely often in bounded intervals.

%s7.1 #&#
\subsection{Exchangeable rewiring process}\label{section:exch rw}
Henceforth, we write $\id$ to denote the identity $\graphsN
\rightarrow\graphsN$ and, for $n\in\mathbb{N}$, we write $\idn$ to
denote the identity $\mathcal{G}_n  \rightarrow
\mathcal{G}_n  $.
Let $\omega$ be an exchangeable measure on $\wireN$ such that
%
%e7.1 #&#
\begin{equation}
\label{eq:regularity omega} \omega\bigl(\{\id\}\bigr)=0 \quad\mbox{and}\quad
 \omega\bigl(\{W\in\wireN
\dvt W|_{[n]}\neq\idn\}\bigr)<\infty \quad\quad\mbox{for every }n\geq2.
\end{equation}
Similar to our definition of $P_{\omega}$ in Section~\ref{section:discrete}, we use $\omega$ to define the transition rates of
continuous-time $\omega$-rewiring chain. Briefly, we assume $\omega(\{
\id\})=0$ because the identity map $\graphsN\rightarrow\graphsN$ is
immaterial for continuous-time processes. The finiteness assumption on
the right of \eqref{eq:regularity omega} ensures that the paths of the
finite restrictions are c\`adl\`ag.

For each $n\in\mathbb{N}$, we write $\omega^{(n)}$ to denote the
restriction of $\omega$ to $ \mathcal{W}_n  $ and define
%
%e7.2 #&#
\begin{equation}
\label{eq:omega-rates} q_{\omega}^{(n)}\bigl(G,G'\bigr):=
\cases{ \omega^{(n)}\bigl(\bigl\{W\in \mathcal{W}_n
\dvt W(G)=G'\bigr\}\bigr),& \quad$G\neq G'\in
\mathcal{G}_n$,
\cr
0,&\quad $G=G'\in\mathcal{G}_n
$. }
\end{equation}
%
%pr7.1 #&#
\begin{prop}\label{prop:finite q-omega}
For each $n\in\mathbb{N}$, $q_{\omega}^{(n)}$ is a finite,
exchangeable conditional measure on $\mathcal{G}_n
$. Moreover, the collection $\{q_{\omega}^{(n)}\}_{n\geq2}$ satisfies
%
%e7.3 #&#
\begin{equation}
\label{eq:consistent rates} q_{\omega}^{(m)}\bigl(G,G'
\bigr)=q_{\omega}^{(n)}\bigl(G^*,\mathbf {R}^{-1}_{m,n}
\bigl(G'\bigr)\bigr),\quad\quad G\neq G'\in\graphsm,
\end{equation}
for all $G^*\in\mathbf{R}^{-1}_{m,n}(G)$, for all $m\leq n$, for
every $n\in\mathbb{N}$, where $\Rmn$ is the restriction map $\mathcal{G}_n  \rightarrow\graphsm$ defined in \eqref
{eq:restriction-graph}.
\end{prop}
\begin{pf}
Finiteness of $q_{\omega}^{(n)}$ follows from \eqref{eq:regularity
omega} since, for every $G\in\mathcal{G}_n  $,
\[
q_{\omega}^{(n)}(G,\mathcal{G}_n )=q_{\omega
}^{(n)}
\bigl(G,\mathcal{G}_n \setminus\{G\}\bigr)=\omega ^{(n)}
\bigl(\bigl\{W\in \mathcal{W}_n \dvt W(G)\neq G\bigr\}\bigr)\leq
\omega^{(n)}\bigl(\{W\neq\idn\}\bigr)<\infty.
\]
Exchangeability of $q_{\omega}^{(n)}$ follows by Proposition~\ref
{prop:exchangeable rewiring chain} and exchangeability of $\omega$.
Consistency of $\{q_{\omega}^{(n)}\}_{n\geq2}$ results from Lipschitz
continuity of rewiring maps (Proposition~\ref{prop:Lipschitz}) and
consistency of the finite-dimensional marginals $\{\omega^{(n)}\}
_{n\in\mathbb{N}}$ associated to $\omega$: for fixed $G\neq G'\in
\graphsm$ and $G^*\in\mathbf{R}^{-1}_{m,n}(G)$,
\begin{eqnarray*}
q_{\omega}^{(n)}\bigl(G^*,\mathbf{R}^{-1}_{m,n}
\bigl(G'\bigr)\bigr)&=&\sum_{G'':G''|_{[m]}=G'}q_{\omega}^{(n)}
\bigl(G^*,G''\bigr)
\\
&=&\sum_{G'':G''|_{[m]}=G'}\omega^{(n)}\bigl(\bigl\{W\in
\mathcal {W}_n \dvt W\bigl(G^*\bigr)=G''\bigr\}
\bigr)
\\
&=&\omega^{(n)}\bigl(\bigl\{W\in \mathcal{W}_n
\dvt W|_{[m]}(G)=G'\bigr\}\bigr)
\\
&=&\omega^{(m)}\bigl(\bigl\{W\in\wirem\dvt W(G)=G'\bigr\}
\bigr)
\\
&=&q_{\omega}^{(m)}\bigl(G,G'\bigr).
\end{eqnarray*}
\upqed\end{pf}

From $\{q_{\omega}^{(n)}\}_{n\in\mathbb{N}}$, we define a collection
of infinitesimal jump rates $\{Q_{\omega}^{(n)}\}_{n\in\mathbb{N}}$ by
%
%e7.4 #&#
\begin{equation}
\label{eq:infinitesimal-omega} Q_{\omega}^{(n)}\bigl(G,G'\bigr):=
\cases{ q_{\omega}^{(n)}\bigl(G,G'\bigr),&\quad
$G'\neq G$,
\cr
-q_{\omega}^{(n)}\bigl(G,
\mathcal{G}_n \setminus\{G\} \bigr),&\quad $G'=G$. }
\end{equation}

%co7.1 #&#
\begin{cor}\label{cor:inf-Q}
The infinitesimal generators $\{Q_{\omega}^{(n)}\}_{n\in\mathbb{N}}$
are exchangeable and consistent and, therefore, define the
infinitesimal jump rates $Q_{\omega}$ of an exchangeable Markov
process on $\graphsN$.
\end{cor}

\begin{pf}
Consistency when $G'\neq G$ was already shown in Proposition~\ref
{prop:finite q-omega}. We must only show that $Q_{\omega}^{(n)}$ is
consistent for $G'=G$.
Fix $n\in\mathbb{N}$ and $G\in\mathcal{G}_n  $.
Then, for any $G^*\in\mathbf{R}^{-1}_{n,n+1}(G)$, we have
\begin{eqnarray*}
Q_{\omega}^{(n+1)}\bigl(G^*,\mathbf{R}^{-1}_{n,n+1}(G)
\bigr)&=&-q_{\omega
}^{(n+1)}\bigl(G^*,\mathcal{G}_{n+1}
\setminus\bigl\{G^*\bigr\}\bigr)+\sum_{G''\in
\mathbf{R}^{-1}_{n,n+1}(G):G''\neq G^*}q_{\omega}^{(n+1)}
\bigl(G^*,G''\bigr)
\\
&=&-q_{\omega}^{(n+1)}\bigl(G^*,\mathcal{G}_{n+1}
\setminus\mathbf {R}^{-1}_{n,n+1}(G)\bigr)
\\
&=&-q_{\omega}^{(n)}\bigl(G,\mathcal{G}_n \setminus\{
G\}\bigr)
\\
&=&Q_{\omega}^{(n)}(G,G).
\end{eqnarray*}
\upqed\end{pf}

In Section~\ref{section:informal description}, we mentioned local and
global discontinuities for graph-valued processes.
In the next two sections, we formally incorporate these discontinuities
into a continuous-time rewiring process: in Section~\ref{section:the
rewiring measure}, we extend the notion of random rewiring from
discrete-time; in Section~\ref{section:local-edge}, we introduce
transitions for which, at the time of a jump, only a single edge in the
network changes. Over time, the local changes can accumulate to cause a
non-trivial change to network structure.

%s7.2 #&#
\subsection{Global jumps: Rewiring}\label{section:the rewiring measure}
In this section, we specialize to the case where $\omega=\Omega
_{\Upsilon}$ for some measure $\Upsilon$ on $\rewiringlimits$ satisfying
%
%e7.5 #&#
\begin{equation}
\label{eq:regularity Upsilon} \Upsilon\bigl(\{\mathbf{I}\}\bigr)=0
\quad\mbox{and}\quad \int
_{\rewiringlimits}\bigl(1-\upsilon_*^{(2)}\bigr)\Upsilon(\mathrm{d}
\upsilon)<\infty,
\end{equation}
where $\mathbf{I}$ is the rewiring limit of $\id\in\wireN$ and
$\upsilon_*^{(n)}:=\upsilon(\idn)$ is the entry of $\upsilon$
corresponding to $\idn$, for each $n\in\mathbb{N}$. For each $n\in
\mathbb{N}$, we write $q_{\Upsilon}^{(n)}$ to denote $q_{\omega
}^{(n)}$ for $\omega=\Omega_{\Upsilon}$, and likewise for the
infinitesimal generator $Q_{\Upsilon}^{(n)}$.

%le7.1 #&#
\begin{lemma}\label{lemma:upsilon-omega equiv}
For $\Upsilon$ satisfying \eqref{eq:regularity Upsilon}, the rewiring
measure $\Omega_{\Upsilon}$ satisfies \eqref{eq:regularity omega}.
\end{lemma}
\begin{pf}
By Theorem~\ref{prop:dissociated}, $\Upsilon(\{\mathbf{I}\})=0$
implies $\Omega_{\Upsilon}(\{\id\})=0$.
We need only show that $\int_{\rewiringlimits}(1-\upsilon
_*^{(2)})\Upsilon(\mathrm{d}\upsilon)<\infty$ implies $\Omega_{\Upsilon
}^{(n)}(\{W\in \mathcal{W}_n  \dvt W\neq\idn\}
)<\infty$ for every $n\geq2$.
For any $\upsilon\in\rewiringlimits$,
\begin{eqnarray*}
\Omega_{\upsilon}\bigl(\{W\in\wireN\dvt W|_{[n]}\neq\idn\}\bigr)&=&
\Omega _{\upsilon} \biggl(\bigcup_{1\leq i<j\leq n}\{W\in
\wireN\dvt W|_{\{i,j\}
}\neq\idij\} \biggr)
\\
&\leq&\sum_{1\leq i<j\leq n}\Omega_{\upsilon}\bigl(\{W\in
\wireN\dvt W|_{\{
i,j\}}\neq\idij\}\bigr)
\\
&=&\sum_{1\leq i<j\leq n}\Omega_{\upsilon}^{(2)}
\bigl(\mathcal {W}_2\setminus\{\idtwo\}\bigr)
\\
&=&\frac{n(n-1)}{2}\bigl(1-\upsilon^{(2)}_*\bigr).
\end{eqnarray*}
Hence, by \eqref{eq:regularity Upsilon},
\[
\Omega_{\Upsilon}\bigl(\{W\in\wireN\dvt W|_{[n]}\neq\idn\}\bigr)\leq
\int_{\rewiringlimits}\frac{n(n-1)}{2}\bigl(1-\upsilon^{(2)}_*
\bigr)\Upsilon (\mathrm{d}\upsilon)<\infty,
\]
for every $n\geq2$.
\end{pf}

%pr7.2 #&#
\begin{prop}
For each $n\in\mathbb{N}$, $q_{\Upsilon}^{(n)}$ is a finite,
exchangeable conditional measure on $\mathcal{G}_n
$. Moreover, $\{q_{\Upsilon}^{(n)}\}_{n\in\mathbb{N}}$ satisfies
\[
q_{\Upsilon}^{(m)}\bigl(G,G'\bigr)=q_{\Upsilon}^{(n)}
\bigl(G^{*},\mathbf {R}^{-1}_{m,n}
\bigl(G'\bigr)\bigr),\quad\quad G\neq G'\in\graphsm, \quad\quad\mbox{for
all }G^*\in\mathbf{R}^{-1}_{m,n}(G).
\]
\end{prop}
\begin{pf}
This follows directly from Propositions \ref{prop:rewiring measure},
\ref{prop:finite q-omega}, and Lemma~\ref{lemma:upsilon-omega equiv}.
\end{pf}

We may, therefore, define an infinitesimal generator for a Markov chain
on $\mathcal{G}_n  $ by
%
%e7.6 #&#
\begin{equation}
\label{eq:infinitesimal-Q} Q_{\Upsilon}^{(n)}\bigl(G,G'\bigr):=
\cases{ q_{\Upsilon}^{(n)}\bigl(G,G'\bigr),&
\quad$G'\neq G$,
\cr
-q_{\Upsilon}^{(n)}\bigl(G,
\mathcal{G}_n \setminus\{ G\}\bigr),&\quad $G'=G$. }
\end{equation}

%th7.1 #&#
\begin{theorem}\label{thm:Upsilon-rewiring process}
For each $\Upsilon$ satisfying \eqref{eq:regularity Upsilon}, there
exists an exchangeable Markov process $\boldsymbol{\Gamma}$ on $\graphsN
$ with finite-dimensional transition rates as in \eqref{eq:infinitesimal-Q}.
\end{theorem}
We call $\boldsymbol{\Gamma}$ in Theorem~\ref{thm:Upsilon-rewiring
process} a {\em rewiring process} directed by $\Upsilon$, or with
rewiring measure~$\Omega_{\Upsilon}$.

%s7.3 #&#
\subsection{Local jumps:~Isolated updating}\label{section:local-edge}
For $i'>j'\geq1$ and $k=0,1$, let $R^k_{i'j'}$ denote the rewiring map
$\wireN\rightarrow\wireN$ that acts by mapping $G\mapsto G':=R^k_{i'j'}(G)$,
%
%e7.7 #&#
\begin{equation}
\label{eq:single edge update} G'_{ij}:=\cases{ G_{ij},&\quad $ij\neq
i'j'$,
\cr
k,&\quad $ij=i'j'$. }
\end{equation}
In words, $R^k_{ij}$ puts an edge between $i$ and $j$ (if $k=1$) or no
edge between $i$ and $j$ (if $k=0$) and keeps every other edge fixed.

For fixed $n\in\mathbb{N}$, let $\mathbf{0}_n\in\mathcal{G}_n  $ denote the {\em empty graph}, that is, the graph
with no edges. We generate a continuous-time process $\boldsymbol{\Gamma
}_0:=\{\Gamma_0(t)\}_{t\geq0}$ on $\mathcal{G}_n
$ as follows. First, we specify a constant $\mathbf{c}_0>0$ and,
independently for each pair $\{i,j\}\in[n]\times[n]$, $i<j$, we
generate i.i.d. random variables $T_{ij}$ from the Exponential
distribution with rate parameter $\mathbf{c}_0$. Given $\{T_{ij}\}
_{i<j}$, we define $\boldsymbol{\Gamma}_0$ by
\[
i\sim_{\Gamma_0(t)} j \quad\Longleftrightarrow\quad T_{ij}<t,
\]
where $i\sim_G j$ denotes an edge between $i$ and $j$ in $G$. Clearly,
$\boldsymbol{\Gamma}_0$ is exchangeable and converges to a unique
stationary distribution $\delta_{\mathbf{1}_n}$, the point mass at
the complete graph $\mathbf{1}_n$. Moreover, the distribution of
$T_*$, the time until absorption in $\mathbf{1}_n$, is simply the law
of the maximum of $n(n-1)/2$ i.i.d. Exponential random variables with
rate parameter $\mathbf{c}_0$.

Conversely, we could consider starting in $\mathbf{1}_n$, the {\em
complete graph}, and generating the above process in reverse. In this
case, we specify $\mathbf{c}_1>0$ and let $\{T_{ij}\}_{i<j}$ be an
i.i.d. collection of Exponential random variables with rate parameter
$\mathbf{c}_1$. We construct $\boldsymbol{\Gamma}_1:=\{\Gamma_1(t)\}
_{t\geq0}$, given $\{T_{ij}\}_{i<j}$, by
 \[
i\sim_{\Gamma_1(t)}j \quad\Longleftrightarrow\quad T_{ij}>t.
\]
For $\mathbf{c}_1=\mathbf{c}_0$, this process evolves exactly as the
{\em complement} of $\boldsymbol{\Gamma}_0$, that is,
\[
\boldsymbol{\Gamma}_1\equalinlaw\bar{\boldsymbol{\Gamma}}_0,
\]
where $\bar{\boldsymbol{\Gamma}}_0:=\{\bar{\Gamma}_0(t)\}_{t\geq0}$
is defined by
\[
i\sim_{\bar{\Gamma}_0(t)} j \quad\Longleftrightarrow\quad i\nsim _{\Gamma_0(t)} j,
\]
for all $i\neq j$ and all $t\geq0$.

It is natural to consider the superposition of $\boldsymbol{\Gamma}_0$
and $\boldsymbol{\Gamma}_1$, which we call a {\em$(\mathbf{c}_0,\mathbf
{c}_1)$-local-edge process}. Let $\mathbf{c}_0,\mathbf{c}_1\geq0$
and let $\delta_{ij}^k$ denote the point mass at the single-edge
update map $R_{ij}^k$. Following Section~\ref{section:exch rw}, we define
%
%e7.8 #&#
\begin{equation}
\label{eq:Omega-c0c1} \Omega_{\mathbf{c}_0,\mathbf{c}_1}:=\mathbf{c}_0\sum
_{i<j}\delta _{ij}^0+
\mathbf{c}_1\sum_{i<j}
\delta_{ij}^1.
\end{equation}

%le7.2 #&#
\begin{lemma}
For $\mathbf{c}_0,\mathbf{c}_1\geq0$, $\Omega_{\mathbf
{c}_0,\mathbf{c}_1}$ defined in \eqref{eq:Omega-c0c1} satisfies
\eqref{eq:regularity omega}.
\end{lemma}
\begin{pf}
Since $\Omega_{\mathbf{c}_0,\mathbf{c}_1}$ only charges the
single-edge update maps, it is clear that it assigns zero mass to the
identity map. Also, for any $n\in\mathbb{N}$, the restriction of
$R_{ij}^k$ to $ \mathcal{W}_n  $ coincides with the
identity $\mathcal{G}_n  \rightarrow\mathcal{G}_n  $ except when $1\leq i<j\leq n$. Hence,
\[
\Omega_{\mathbf{c}_0,\mathbf{c}_1}^{(n)}\bigl(\{W\in \mathcal {W}_n \dvt W
\neq\idn\}\bigr)=\frac{n(n-1)}{2}(\mathbf {c}_0+
\mathbf{c}_1)<\infty,
\]
for every $n\geq2$.
\end{pf}

%co7.2 #&#
\begin{cor}
For any $\mathbf{c}_0,\mathbf{c}_1\geq0$, there exists an
exchangeable Markov process on $\graphsN$ with jump rates given by
$\Omega_{\mathbf{c}_0,\mathbf{c}_1}$.
\end{cor}
\begin{pf}
For every $n\in\mathbb{N}$, the total jump rate out of any $G\in
\mathcal{G}_n  $ can be no larger than
\[
\frac{n(n-1)}{2}(\mathbf{c}_0\vee\mathbf{c}_1)<\infty,
\]
and so the finite-dimensional hold times are almost surely positive and
the process on $\mathcal{G}_n  $ has c\`adl\`ag
sample paths. The Markov property and exchangeability follow by
independence of the Exponential hold times $\{T_{ij}\}_{1\leq i<j\leq
n}$ and Corollary~\ref{cor:inf-Q}. Consistency is apparent by the
construction from independent Poisson point processes. This completes
the proof.
\end{pf}

%de7.1 #&#
\begin{defn}
For any measure $\Upsilon$ satisfying \eqref{eq:regularity Upsilon},
$\mathbf{c}_0,\mathbf{c}_1\geq0$, we call a rewiring process with
jump measure $\omega=\Omega_{\Upsilon}+\Omega_{\mathbf
{c}_0,\mathbf{c}_1}$ an {\em$(\Upsilon,\mathbf{c}_0,\mathbf
{c}_1)$-rewiring process}.
\end{defn}

From our discussion in this section, the $(\Upsilon,\mathbf
{c}_0,\mathbf{c}_1)$-rewiring process exists for any choice of
$\Upsilon$ satisfying \eqref{eq:regularity Upsilon} and $\mathbf
{c}_0,\mathbf{c}_1\geq0$. Individually, $\Omega_{\Upsilon}$ and
$\Omega_{\mathbf{c}_0,\mathbf{c}_1}$ satisfy \eqref{eq:regularity
omega} and, thus, so does $\omega:=\Omega_{\Upsilon}+\Omega
_{\mathbf{c}_0,\mathbf{c}_1}$. Furthermore, the family of $(\Upsilon
,\mathbf{c}_0,\mathbf{c}_1)$-rewiring processes is Markovian,
exchangeable, and consistent.

%s8 #&#
\section{Simulating rewiring processes}\label{section:Poissonian structure}
We can construct an $(\Upsilon,\mathbf{c}_0,\mathbf{c}_1)$-rewiring
process from a Poisson point process.
For $\omega:=\Omega_{\Upsilon}+\Omega_{\mathbf{c}_0,\mathbf
{c}_1}$, where $\Upsilon$ satisfies \eqref{eq:regularity Upsilon} and
$\mathbf{c}_0,\mathbf{c}_1\geq0$, let $\mathbf{W}:=\{(t,W_t)\}
\subset\mathbb{R}^+\times\wireN$ be a Poisson point process with
intensity $\mathrm{d}t\otimes\omega$. To begin, we take $\Gamma_0$ to be an
exchangeable random graph and, for each $n\in\mathbb{N}$, we define
$\boldsymbol{\Gamma}^{[n]}:=(\Gamma^{[n]}_t)_{t\geq0}$ on $\mathcal{G}_n  $ by $\Gamma^{[n]}_0:=\Gamma_0|_{[n]}$ and
\begin{itemize}
\item if $t>0$ is an atom time of $\mathbf{W}$ such that
$W^{[n]}_t:=W_t|_{[n]}\neq\idn$, then we put $\Gamma
^{[n]}_{t}:=W^{[n]}_t(\Gamma^{[n]}_{t-})$;
\item otherwise, we put $\Gamma^{[n]}_t=\Gamma^{[n]}_{t-}$.
\end{itemize}

%pr8.1 #&#
\begin{prop}\label{prop:thinned PPP}
For each $n\in\mathbb{N}$, $\boldsymbol{\Gamma}^{[n]}$ is a Markov
chain on $\mathcal{G}_n  $ with infinitesimal jump
rates $Q_{\omega}^{(n)}$ in \eqref{eq:infinitesimal-omega}.
\end{prop}
\begin{pf}
We can define $\mathbf{W}^{[n]}:=\{(t,W_t^{[n]})\}\subset\mathbb
{R}^{+}\times \mathcal{W}_n  $ from $\mathbf{W}$
by removing any atom times for which $W_t^{[n]}:=W_t|_{[n]}=\idn$, and
otherwise putting $W_t^{[n]}:=W_t|_{[n]}$. By the thinning property of
Poisson point processes, $\mathbf{W}^{[n]}$ is a Poisson point process
with intensity $\mathrm{d}t\otimes\omega_n$, where
\[
\omega_n(\cdot):=\omega^{(n)}\bigl(\cdot\setminus\{\idn\}
\bigr).
\]
Given $\Gamma_t^{[n]}=G$, the jump rate to state $G'\neq G$ is
\[
\omega_n\bigl(\bigl\{W\in \mathcal{W}_n
\dvt W(G)=G'\bigr\} \bigr)=Q_{\omega}^{(n)}
\bigl(G,G'\bigr),
\]
and the conclusion follows.
\end{pf}

%th8.1 #&#
\begin{theorem}\label{thm:existence rewiring}
For any $\omega$ satisfying \eqref{eq:regularity omega}, the $\omega
$-rewiring process on $\graphsN$ exists and can be constructed from a
Poisson point process with intensity $\mathrm{d}t\otimes\omega$ as above.
\end{theorem}

\begin{pf}
Let $\mathbf{W}$ be a Poisson point process with intensity $\mathrm{d}t\otimes
\omega$ and construct $\{\boldsymbol{\Gamma}^{[n]}\}_{n\in\mathbb{N}}$
from the thinned processes $\{\mathbf{W}^{[n]}\}_{n\in\mathbb{N}}$
determined by $\mathbf{W}$. By Proposition~\ref{prop:thinned PPP},
each $\boldsymbol{\Gamma}^{[n]}$ is an exchangeable Markov chain governed
by $Q_{\omega}^{(n)}$. Moreover, $\{\boldsymbol{\Gamma}^{[n]}\}_{n\in
\mathbb{N}}$ is compatible by construction, that is, $\Gamma
^{[m]}_t=\Rmn\Gamma^{[n]}_t$ for all $t\geq0$, for all $m\leq n$;
hence, $\{\boldsymbol{\Gamma}^{[n]}\}_{n\in\mathbb{N}}$ defines a
process $\boldsymbol{\Gamma}$ on $\graphsN$. As we have shown
previously, the infinitesimal rates given by $\{Q_{\omega}^{(n)}\}
_{n\in\mathbb{N}}$ are consistent and exchangeable; hence, $\boldsymbol
{\Gamma}$ has infinitestimal generator $Q_{\omega}$ and is an $\omega
$-rewiring process.
\end{pf}

%s8.1 #&#
\subsection{The Feller property}\label{section:Feller}
Any Markov process $\boldsymbol{\Gamma}$ on $\graphsN$ is characterized
by its semigroup $(\mathbf{P}_t)_{t\geq0}$, defined as an operator on
the space of continuous, bounded functions $h\dvtx \graphsN\rightarrow
\mathbb{R}$ by
%
%e8.1 #&#
\begin{equation}
\label{eq:semigroup} \mathbf{P}_th(G):=\mathbb{E}_{G}h(
\Gamma_t),\quad\quad G\in\graphsN,
\end{equation}
where $\mathbb{E}_G$ denotes the expectation operator with respect to
the initial distribution $\delta_G(\cdot)$, the point mass at $G$. We
say $\boldsymbol{\Gamma}$ has the {\em Feller property} if, for all
bounded, continuous functions $h\dvtx \graphsN\rightarrow\mathbb{R}$, its
semigroup satisfies
\begin{itemize}
\item $\mathbf{P}_th(G)\rightarrow h(G)$ as $t\downarrow0$ for all
$G\in\graphsN$, and
\item$G\mapsto\mathbf{P}_th(G)$ is continuous for all $t\geq0$.
\end{itemize}

%th8.2 #&#
\begin{theorem}
The semigroup $(\mathbf{P}^{\omega}_t)_{t\geq0}$ of any $\omega
$-rewiring process enjoys the Feller property.
\end{theorem}
\begin{pf}
To show the first point in the Feller property, we let $G\in\graphsN$
and $\boldsymbol{\Gamma}:=(\Gamma_t)_{t\geq0}$ be an $\omega$-rewiring
process with initial state $\Gamma_0=G$ and directing measure $\omega
$ satisfying \eqref{eq:regularity omega}. We define
\[
\mathcal{F}:=\bigl\{h\dvtx \graphsN\rightarrow\mathbb{R}\mid\mbox{there exists }n\in
\mathbb{N}\mbox{ such that }G|_{[n]}=G'|_{[n]}
\Rightarrow h(G)=h\bigl(G'\bigr)\bigr\}.
\]
By \eqref{eq:regularity omega} and finiteness of $\mathcal{G}_n  $, $\Gamma^{[n]}_t\rightarrow G|_{[n]}$ in
probability as $t\downarrow0$, for every $n\in\mathbb{N}$. Thus, for
any $h\in\mathcal{F}$, let $N\in\mathbb{N}$ be such that
\[
d\bigl(G,G'\bigr)\leq1/N \quad\Longrightarrow\quad h(G)=h
\bigl(G'\bigr).
\]
Then $\Gamma_t^{[N]}\rightarrow G|_{[N]}$ in probability as
$t\downarrow0$ and, therefore, $\mathbf{P}_th(G)\rightarrow h(G)$ by
the Bounded Convergence theorem. Right-continuity at zero for all
bounded, continuous $h\dvtx \graphsN\rightarrow\mathbb{R}$ follows by the
Stone--Weierstrass theorem.

For the second point, let $G,G'\in\graphsN$ have $d(G,G')\leq1/n$
for some $n\in\mathbb{N}$ and construct $\boldsymbol{\Gamma}$ and
$\boldsymbol{\Gamma}'$ from the same Poisson point process $\mathbf{W}$
but with initial states $\Gamma_0=G$ and $\Gamma'_0=G'$. By Lipschitz
continuity of the rewiring maps (Proposition~\ref{prop:Lipschitz}),
$\boldsymbol{\Gamma}$ and $\boldsymbol{\Gamma}'$ can never be more than
distance $1/n$ apart, for all $t\geq0$. Continuity of $\mathbf
{P}_t^{\omega}$, for each $t\geq0$, follows.
\end{pf}

By the Feller property, any $\omega$-rewiring process has a c\`adl\`ag
version and its jumps are characterized by an infinitesimal generator.
In Section~\ref{section:continuous}, we described the infinitesimal
generator through its finite restrictions. Ethier and Kurtz \cite
{EthierKurtz} give an extensive treatment of the general theory of
Feller processes.

%s9 #&#
\section{Concluding remarks}\label{section:concluding remarks}

We have presented a family of time-varying network models that is
Markovian, exchangeable, and consistent, natural statistical properties
that impose structure without introducing logical pitfalls.
External to statistics, exchangeable models are flawed: they produce
{\em dense} graphs when conventional wisdom suggests real-world
networks are {\em sparse}.
The Erd\H{o}s--R\'enyi model's storied history cautions against dismay.
Though it replicates little real-world network structure, the Erd\H
{o}s--R\'enyi model has produced a deluge of insight for graph-theoretic
structures and is a paragon of the utility of the probabilistic method
\cite{AlonSpencerBook}.
While our discussion is specific to exchangeable processes, the general
descriptions in Sections~\ref{section:discrete} and \ref
{section:continuous} can be used to construct processes that are not
exchangeable, and possibly even sparse.

The most immediate impact of the rewiring process may be for analyzing
information spread on dynamic networks.
Under the heading of {\em Finite Markov Information Exchange} (FMIE)
processes, Aldous \cite{AldousFMIE2013} recently surveyed interacting
particle systems models for social network dynamics.
Informally, FMIE processes model a random spread of information on a network.
Some of the easiest to describe FMIE processes coincide with well-known
interacting particle systems, such as the Voter and Contact processes;
others mimic certain social behaviors, for example, {\em Fashionista}
and {\em Compulsive Gambler}.

Simulation is a valuable practical tool for developing intuition about
intractable problems.
Aldous's expository account contains some hard open problems for
time-invariant networks.
Considering these same questions on dynamic networks seems an even
greater challenge.
Despite these barriers, policymakers and scientists alike desire to
understand how trends, epidemics, and other information spread on networks.
The Poisson point process construction in Section~\ref{section:Poissonian structure} could be fruitful for deriving practical
answers to these problems.

%s10 #&#
\section{Technical proofs}\label{section:proof}
In this section, we prove some technical results from our previous discussion.

%s10.1 #&#
\subsection{Proof of Lemma \texorpdfstring{\protect\ref{lemma:compact
rewiring}}{6.2}}\label{section:compact rewiring}
We now show that $(\rewiringlimits,\rho)$ is a compact metric space.
Recall that $\rewiringlimits$ is equipped with the metric
\[
\rho\bigl(\upsilon,\upsilon'\bigr)=\sum
_{n\in\mathbb{N}}2^{-n}\sum_{V\in
\mathcal{W}_n  }\bigl|
\upsilon_V-\upsilon'_V\bigr|, \quad\quad\upsilon,
\upsilon'\in\rewiringlimits.
\]
Since $[0,1]^{\mathcal{W}^*}$ is compact in this metric, it suffices
to show that $\rewiringlimits$ is a closed subset of $[0,1]^{\mathcal{W}^*}$.

By Lemma~\ref{lemma:rewiring structure}, every $\upsilon\in
\rewiringlimits$ satisfies
\[
\upsilon(W)=\sum_{W^*\in\mathcal{W}_{n+1}:W^*|_{[n]}=W}\upsilon \bigl(W^*\bigr)
\quad\quad\mbox{for every }W\in \mathcal{W}_n
\]
and
\[
\sum_{W\in \mathcal{W}_n  }\upsilon(W)=1,
\]
for all $n\in\mathbb{N}$.
Then, for any $x\in[0,1]^{\mathcal{W}^*}\setminus\rewiringlimits$,
there must be some $N\in\mathbb{N}$ for which
\[
\varepsilon_x:=\sum_{W\in\mathcal{W}_N}\biggl\llvert
x^{(N)}(W)-\sum_{W^*|_{[N]}=W}x^{(N+1)}
\bigl(W^*\bigr)\biggr\rrvert >0.
\]
For any $\delta>0$, let $B(x,\delta):=\{x'\in[0,1]^{\mathcal
{W}^*}\dvt  \rho(x,x')<\delta\}$ denote the $\delta$-ball around $x$.
Now, take any $x'\in B(x,2^{-N-2}\varepsilon_x)$. By this assumption,
$\rho(x,x')\leq2^{-N-2}\varepsilon_x$ and so
\begin{eqnarray*}
&&2^{-N}\sum_{W\in\mathcal{W}_N}\bigl|x^{(N)}(W)-x'^{(N)}(W)\bigr|+2^{-N-1}
\sum_{W^*\in\mathcal{W}_{N+1}}\bigl|x^{(N+1)}\bigl(W^*
\bigr)-x'^{(N+1)}\bigl(W^*\bigr)\bigr|\\
&&\quad\leq 2^{-N-2}
\varepsilon_x;
\end{eqnarray*}
whence,
\begin{eqnarray*}
\sum_{W\in\mathcal{W}_N}\bigl|x^{(N)}(W)-x'^{(N)}(W)\bigr|
&\leq&\frac
{1}{4}\varepsilon_x \quad\mbox{and}
\\
\sum_{W^*\in\mathcal{W}_{N+1}}\bigl|x^{(N+1)}\bigl(W^*
\bigr)-x'^{(N+1)}\bigl(W^*\bigr)\bigr|&\leq& \frac{1}{2}
\varepsilon_x.
\end{eqnarray*}
By the triangle inequality, we have
\begin{eqnarray*}
\varepsilon_x&=&\sum_{W\in\mathcal{W}_N}\biggl
\llvert x^{(N)}(W)-\sum_{W^*:W^*|_{[N]}}x^{(N+1)}
\bigl(W^*\bigr)\biggr\rrvert
\\
&\leq&\sum_{W\in\mathcal{W}_N}\bigl\llvert x^{(N)}(W)-x'^{(N)}(W)
\bigr\rrvert +\sum_{W\in\mathcal{W}_N}\biggl\llvert \sum
_{W^*:W^*|_{[N]}=W}\bigl(x^{(N+1)}\bigl(W^*\bigr)-x'^{(N+1)}
\bigl(W^*\bigr)\bigr)\biggr\rrvert
\\
&&{} +\sum_{W\in\mathcal{W}_N}\biggl\llvert x'^{(N)}(W)-
\sum_{W^*:W^*|_{[N]}=W}x'^{(N+1)}\bigl(W^*
\bigr)\biggr\rrvert
\\
&\leq&\varepsilon_x/4+\sum_{W\in\mathcal{W}_N}\sum
_{W^*:W^*|_{[N]}=W}\bigl\llvert x^{(N+1)}\bigl(W^*
\bigr)-x'^{(N+1)}\bigl(W^*\bigr)\bigr\rrvert
\\
&&{}+\sum
_{W\in\mathcal{W}_N}\biggl\llvert x'^{(N)}(W)-\sum
_{W^*:W^*|_{[N]}=W}x'^{(N+1)}\bigl(W^*\bigr)
\biggr\rrvert
\\
&\leq&\varepsilon_x/4+\varepsilon_x/2+\sum
_{W\in\mathcal
{W}_N}\biggl\llvert x'^{(N)}(W)-\sum
_{W^*:W^*|_{[N]}=W}x'^{(N+1)}\bigl(W^*\bigr)
\biggr\rrvert .
\end{eqnarray*}
Therefore,
\[
\sum_{W\in\mathcal{W}_N}\biggl\llvert x'^{(N)}(W)-
\sum_{W^*:W^*|_{[N]}=W}x'^{(N+1)}\bigl(W^*
\bigr)\biggr\rrvert \geq\varepsilon_x/4>0,
\]
which implies $x'\in[0,1]^{\mathcal{W}^*}\setminus\rewiringlimits$,
meaning $[0,1]^{\mathcal{W}^*}\setminus\rewiringlimits$ is open and
$\rewiringlimits$ is closed. Since $[0,1]^{\mathcal{W}^*}$ is
compact, so is $\rewiringlimits$. This completes the proof.

%s10.2 #&#
\subsection{Proof of Theorem \texorpdfstring{\protect\ref{prop:dissociated}}{6.1}}
Assume that $W$ is an exchangeable and dissociated rewiring map. By the
Aldous--Hoover theorem, we can assume $W$ is constructed from a
measurable function $f\dvtx [0,1]^4\rightarrow\{0,1\}\times\{0,1\}$ for
which (i) $f(a,\cdot,\cdot,\cdot)=f(a',\cdot,\cdot,\cdot)$ and
(ii) $f(\cdot,b,c,\cdot)=f(\cdot,c,b,\cdot)$. More precisely, we
assume $W_{ij}=f(\alpha,\eta_i,\eta_j,\lambda_{\{i,j\}})$, for each
$i,j\geq1$, where $\{\alpha;(\eta_i)_{i\geq1};(\lambda_{\{i,j\}
})_{i<j}\}$ are i.i.d. Uniform random variables on $[0,1]$.
Conditional on $\alpha=a$, we define the quantity
\[
t_a(V,W):=P\{W|_{[m]}=V\mid\alpha=a\},\quad\quad V\in\wirem, m\in
\mathbb{N},
\]
which, by the fact that $W$ is dissociated, is independent of $a$;
hence, we define the non-random quantity
\[
t(V,W):=E\bigl(\mathbf{1}\{W|_{[m]}=V\}\mid\alpha\bigr)=P
\{W|_{[m]}=V\}.
\]
Recall, from Section~\ref{section:rewiring limits}, the definition
\[
t(V,W|_{[n]}):=\frac{\ind(V,W|_{[n]})}{n^{\downarrow m}}:=\frac
{1}{n^{\downarrow m}}\sum
_{\mathrm{injections }\psi:[m]\rightarrow
[n]}\mathbf{1}\bigl\{W|_{[n]}^{\psi}=V\bigr
\},\quad\quad n\in\mathbb{N}.
\]
For every $n\geq1$, we also define
\[
M_{k,n}:=\frac{1}{n^{\downarrow m}}\sum_{\mathrm{injections }\psi
:[m]\rightarrow[n]}E
\bigl(\mathbf{1}\bigl\{W|_{[n]}^{\psi}=V\bigr\} \mid
W|_{[k]}\bigr),\quad\quad k=0,1,\ldots,n.
\]
In particular, for every $n\in\mathbb{N}$, we have
\[
M_{0,n}=\frac{1}{n^{\downarrow m}}\sum_{\mathrm{injections }\psi
:[m]\rightarrow[n]}E
\bigl(\mathbf{1}\bigl\{W|_{[n]}^{\psi}=V\bigr\} \mid
W|_{[0]}\bigr)=t(V,W)
\]
and
\[
M_{n,n}=\frac{1}{n^{\downarrow m}}\sum_{\mathrm{injections }\psi
:[m]\rightarrow[n]}E
\bigl(\mathbf{1}\bigl\{W|^{\psi}_{[n]}=V\bigr\} \mid W|_{[n]}
\bigr)=t(V,W|_{[n]}).
\]
We wish to show that $t(V,W|_{[n]})\rightarrow t(V,W)$ almost surely,
for every $V\in\wirem$, $m\in\mathbb{N}$. To do this, we first show
that $(M_{0,n},M_{1,n},\ldots,M_{n,n})$ is a martingale with respect
to its natural filtration, for every $n\in\mathbb{N}$. We can then
appeal to Azuma's inequality and the Borel--Cantelli lemma to show that
$M_{n,n}\rightarrow t(V,W)$ as $n\rightarrow\infty$.

Note that
\[
M_{k,n}=\frac{1}{n^{\downarrow m}}\sum_{\mathrm{injections }\psi
:[m]\rightarrow[n]}\sum
_{w\in \mathcal{W}_n
}E\bigl(\mathbf{1}\{W|_{[n]}=w\} \mid
{W}|_{[k]}\bigr)\mathbf{1}\bigl\{w^{\psi}=V\bigr\}
\]
and
\[
E(M_{k+1,n}\mid M_{k,n})=E\bigl(E(M_{k+1,n}\mid
M_{k,n},{W}|_{[k]}) \mid M_{k,n}\bigr).
\]
On the inside, we have
\begin{eqnarray*}
&&E(M_{k+1,n}\mid M_{k,n},W|_{[k]})
\\
&&\quad=E \biggl(
\frac{1}{n^{\downarrow
m}}\sum_{\mathrm{injections\ }\psi:[m]\rightarrow[n]}\sum
_{w\in
\mathcal{W}_n  }\mathbf{1}\bigl\{w^{\psi}=V\bigr\}E\bigl(\mathbf
{1}\bigl\{W|_{[n]}^{\psi}=w\bigr\}\mid W|_{[k]}\bigr) \Bigm|
M_{k,n},W|_{[k]} \biggr)
\\
&&\quad=\frac{1}{n^{\downarrow m}}\sum_{\mathrm{injections\ }\psi
:[m]\rightarrow[n]}\sum
_{w\in \mathcal{W}_n
}\mathbf{1}\bigl\{w^{\psi}=V\bigr\}E\bigl(
\mathbf{1}\bigl\{W|_{[n]}^{\psi}=w\bigr\} \mid W|_{[k]}
\bigr);
\end{eqnarray*}
whence,
\begin{eqnarray*}
&&E\bigl(E(M_{k+1,n}\mid M_{k,n},W|_{[k]})\mid
M_{k,n}\bigr)
\\
&&\quad=E \biggl(\frac
{1}{n^{\downarrow m}}\sum
_{\mathrm{injections\ }\psi:[m]\rightarrow
[n]}\sum_{w\in \mathcal{W}_n  }\mathbf{1}\bigl\{
w^{\psi}=V\bigr\}E\bigl(\mathbf{1}\bigl\{W|_{[n]}^{\psi}=w
\bigr\}\mid W|_{[k]}\bigr) \Bigm| M_{k,n} \biggr)
\\
&&\quad=\frac{1}{n^{\downarrow m}}\sum_{\mathrm{injections\ }\psi
:[m]\rightarrow[n]}\sum
_{w\in \mathcal{W}_n
}\mathbf{1}\bigl\{w^{\psi}=V\bigr\}E \bigl(E
\bigl(\mathbf{1}\{W|_{[n]}=w\}\mid W|_{[k]}\bigr) \mid
M_{k,n} \bigr)
\\
&&\quad=\frac{1}{n^{\downarrow m}}\sum_{\mathrm{injections\ }\psi
:[m]\rightarrow[n]}\sum
_{w\in \mathcal{W}_n
}\mathbf{1}\bigl\{w^{\psi}=V\bigr\}E\bigl(
\mathbf{1}\{W|_{[n]}=w\} \mid M_{k,n}\bigr)
\\
&&\quad=M_{k,n}.
\end{eqnarray*}
Therefore, $(M_{k,n})_{k=0,1,\ldots,n}$ is a martingale for every
$n\in\mathbb{N}$. Furthermore, for every $k=0,1,\ldots,n-1$,
\begin{eqnarray*}
&&|M_{k+1,n}-M_{k,n}|
\\
&&\quad=\frac{1}{n^{\downarrow m}}\biggl\llvert \sum
_{\mathrm{injections\ }\psi:[m]\rightarrow[n]}E\bigl(\mathbf{1}\bigl\{W|_{[n]}^{\psi
}=V
\bigr\}\mid W|_{[k+1]}\bigr)-E\bigl(\mathbf{1}\bigl\{W|_{[n]}^{\psi}=V
\bigr\}\mid W|_{[k]}\bigr)\biggr\rrvert
\\
&&\quad\leq\frac{1}{n^{\downarrow m}}\sum_{\mathrm{injections\ }\psi
:[m]\rightarrow[n]}\bigl|E\bigl(
\mathbf{1}\bigl\{W|_{[n]}^{\psi}=V\bigr\}\mid W|_{[k+1]}
\bigr)-E\bigl(\mathbf{1}\bigl\{W|_{[n]}^{\psi}=V\bigr\}\mid
W|_{[k]}\bigr)\bigr|
\\
&&\quad\leq m(n-1)^{\downarrow(m-1)}/n^{\downarrow m}
\\
&&\quad\leq m/n,
\end{eqnarray*}
since $E(\mathbf{1}\{W|_{[n]}^{\psi}=V\}\mid W|_{[k+1]})-E(\mathbf
{1}\{W|_{[n]}^{\psi}=V\}\mid W|_{[k]})=0$ whenever $\psi$ does not
map an element to $k+1$. The conditions for Azuma's martingale
inequality are thus satisfied and we have, for every $\varepsilon>0$,
\[
P\bigl\{|M_{n,n}-M_{0,n}|>\varepsilon\bigr\}\leq2\exp \biggl\{-
\frac
{\varepsilon^{2}n}{2m^2} \biggr\} \quad\quad\mbox{for every }n\in\mathbb{N}.
\]
Thus,
\[
\sum_{n=1}^{\infty}P\bigl\{\bigl|M_{n,n}-t(V,W)\bigr|>
\varepsilon\bigr\}\leq2\sum_{n=1}^{\infty}\exp
\biggl\{-\frac{\varepsilon^2n}{2m^2} \biggr\} <\infty,
\]
and we conclude, by the Borel--Cantelli lemma, that
\begin{eqnarray*}
&&\limsup_{n\rightarrow\infty} \bigl\{ \bigl|t(V,W|_{[n]})-t(V,W)\bigr|>
\varepsilon \bigr\}
\\
&&\quad= \bigl\{ \bigl|t(V,W|_{[n]})-t(V,W)\bigr|>\varepsilon\mbox{ for
infinitely many }n\in \mathbb{N} \bigr\}
\end{eqnarray*}
has probability zero. It follows that $\lim_{n\rightarrow\infty
}t(V,W|_{[n]})=t(V,W)$ exists with probability one for every $V\in
\bigcup_{m\in\mathbb{N}}\wirem$. Therefore, with probability one,
the rewiring limit $(t(V,W))_{V\in\mathcal{W}^*}$ exists. We have
already shown, by the assumption that $W$ is dissociated, that $t(V,W)$
is non-random for every $V\in\bigcup_{m\in\mathbb{N}}\wirem$;
hence, the limit $(t(V,W))_{V\in\mathcal{W}^*}$ is non-random. This
completes the proof.
%\begin{appendix}
%\section{}
%\end{appendix}

\section*{Acknowledgement}
This work is partially supported by NSF Grant DMS-1308899 and NSA Grant
H98230-13-1-0299.

%\begin{supplement}%[id=suppA]
%\slink[doi]{10.3150/00-BEJXXXXSUPP} %[doi,text={...}] - jei reikia
%suskaldyti doi
%\sfilename{BEJ000\_supp.pdf}
%\end{supplement}

% imsref loaded by arune.pranskunaite, 2014-04-24 13:08:52
%
% imsref loaded by arune.pranskunaite, 2014-05-05 09:15:38

\printhistory


\begin{thebibliography}{31}

%%% bbsrt2.pl, ver. 2.5.2, 2014.04.15

%b2 ###bbsrt2
\bibitem{AldousExchArrays}
\begin{barticle}[mr]
\bauthor{\bsnm{Aldous},~\bfnm{David~J.}\binits{D.J.}}
(\byear{1981}).
\btitle{Representations for partially exchangeable arrays of random variables}.
\bjournal{J. Multivariate Anal.}
\bvolume{11}
\bpages{581--598}.
\bid{doi={10.1016/0047-259X(81)90099-3}, issn={0047-259X}, mr={0637937}}
\end{barticle}
\bptok{imsref}\endbibitem

%b3 ###bbsrt2
\bibitem{AldousExchangeability}
\begin{bincollection}[mr]
\bauthor{\bsnm{Aldous},~\bfnm{David~J.}\binits{D.J.}}
(\byear{1985}).
\btitle{Exchangeability and related topics}.
In \bbooktitle{\'{E}cole D'\'et\'e de Probabilit\'es de
{S}aint-{F}lour, {XIII}---1983}.
\bseries{Lecture Notes in Math.}
\bvolume{1117}
\bpages{1--198}.
\bpublisher{Springer, Berlin}.
\bid{doi={10.1007/BFb0099421}, mr={0883646}}
\end{bincollection}
\bptok{imsref}\endbibitem


%b1 ###bbsrt2
\bibitem{Aldous1996cladograms}
\begin{bincollection}[mr]
\bauthor{\bsnm{Aldous},~\bfnm{David}\binits{D.}}
(\byear{1996}).
\btitle{Probability distributions on cladograms}.
In \bbooktitle{Random Discrete Structures ({M}inneapolis, {MN}, 1993)}.
\bseries{IMA Vol. Math. Appl.}
\bvolume{76}
\bpages{1--18}.
\bpublisher{Springer, New York}.
\bid{doi={10.1007/978-1-4612-0719-1_1}, mr={1395604}}
\end{bincollection}
\bptok{imsref}\endbibitem

%b4 ###bbsrt2
\bibitem{AldousFMIE2013}
\begin{bmisc}[auto:STB|2014/02/12|14:17:21]
\bauthor{\bsnm{Aldous},~\bfnm{D.~J.}\binits{D.J.}}
(\byear{2012}).
\bhowpublished{Interacting particle systems as stochastic social dynamics.
Preprint}.
\end{bmisc}
\bptok{imsref}\endbibitem

%b5 ###bbsrt2
\bibitem{AlonSpencerBook}
\begin{bbook}[mr]
\bauthor{\bsnm{Alon},~\bfnm{Noga}\binits{N.}} \AND
\bauthor{\bsnm{Spencer},~\bfnm{Joel~H.}\binits{J.H.}}
(\byear{2000}).
\btitle{The Probabilistic Method},
\bedition{2nd} ed.
\bseries{Wiley-Interscience Series in Discrete Mathematics and Optimization}.
\blocation{New York}:
\bpublisher{Wiley-Interscience}.
\bnote{With an appendix on the life and work of Paul Erd{\H{o}}s}.
\bid{doi={10.1002/0471722154}, mr={1885388}}
\end{bbook}
\bptok{imsref}\endbibitem

%b6 ###bbsrt2
\bibitem{Linked}
\begin{bbook}[auto:STB|2014/02/12|14:17:21]
\bauthor{\bsnm{Barab{\'a}si},~\bfnm{A.-L.}\binits{A.-L.}}
(\byear{2003}).
\btitle{Linked: How Everything Is Connected to Everything Else
and What It Means for Business, Science, and Everyday Life}.
\blocation{New York}:
\bpublisher{Plume}.
\end{bbook}
\bptok{imsref}\endbibitem

%b7 ###bbsrt2
\bibitem{BarabasiAlbert1999}
\begin{barticle}[mr]
\bauthor{\bsnm{Barab{\'a}si},~\bfnm{Albert-L{\'a}szl{\'o}}\binits
{A.-L.}} \AND
\bauthor{\bsnm{Albert},~\bfnm{R{\'e}ka}\binits{R.}}
(\byear{1999}).
\btitle{Emergence of scaling in random networks}.
\bjournal{Science}
\bvolume{286}
\bpages{509--512}.
\bid{doi={10.1126/science.286.5439.509}, issn={0036-8075}, mr={2091634}}
\end{barticle}
\bptok{imsref}\endbibitem

%b8 ###bbsrt2
\bibitem{BurkeRosenblatt1958}
\begin{barticle}[mr]
\bauthor{\bsnm{Burke},~\bfnm{C.~J.}\binits{C.J.}} \AND
\bauthor{\bsnm{Rosenblatt},~\bfnm{M.}\binits{M.}}
(\byear{1958}).
\btitle{A {M}arkovian function of a {M}arkov chain}.
\bjournal{Ann. Math. Statist.}
\bvolume{29}
\bpages{1112--1122}.
\bid{issn={0003-4851}, mr={0101557}}
\end{barticle}
\bptok{imsref}\endbibitem

%b9 ###bbsrt2
\bibitem{ChungLubook}
\begin{bbook}[mr]
\bauthor{\bsnm{Chung},~\bfnm{Fan}\binits{F.}} \AND
\bauthor{\bsnm{Lu},~\bfnm{Linyuan}\binits{L.}}
(\byear{2006}).
\btitle{Complex Graphs and Networks}.
\bseries{CBMS Regional Conference Series in Mathematics}
\bvolume{107}.
\blocation{Providence, RI}:
\bpublisher{Published for the Conference Board of the Mathematical
Sciences, Washington, DC; by the Amer. Math. Soc.}
\bid{mr={2248695}}
\end{bbook}
\bptok{imsref}\endbibitem

%b10 ###bbsrt2
\bibitem{ClarkeCyber}
\begin{bbook}[auto:STB|2014/02/12|14:17:21]
\bauthor{\bsnm{Clarke},~\bfnm{R.}\binits{R.}}
(\byear{2010}).
\btitle{Cyber War: The Next Threat to National Security and What to Do
About It}.
\blocation{New York}:
\bpublisher{HarperCollins}.
\end{bbook}
\bptok{imsref}\endbibitem

%b11 ###bbsrt2
\bibitem{DorogovtsevMendes2003}
\begin{bbook}[mr]
\bauthor{\bsnm{Dorogovtsev},~\bfnm{S.~N.}\binits{S.N.}} \AND
\bauthor{\bsnm{Mendes},~\bfnm{J.~F.~F.}\binits{J.F.F.}}
(\byear{2003}).
\btitle{Evolution of Networks: From Biological Nets to the Internet and WWW}.
\blocation{Oxford}:
\bpublisher{Oxford Univ. Press}.
\bid{doi={10.1093/acprof:oso/9780198515906.001.0001}, mr={1993912}}
\end{bbook}
\bptok{imsref}\endbibitem

%b12 ###bbsrt2
\bibitem{DurrettRandomGraphs}
\begin{bbook}[mr]
\bauthor{\bsnm{Durrett},~\bfnm{Rick}\binits{R.}}
(\byear{2007}).
\btitle{Random Graph Dynamics}.
\bseries{Cambridge Series in Statistical and Probabilistic Mathematics}.
\blocation{Cambridge}:
\bpublisher{Cambridge Univ. Press}.
\bid{mr={2271734}}
\end{bbook}
\bptok{imsref}\endbibitem

%b13 ###bbsrt2
\bibitem{ErdosRenyi1959}
\begin{barticle}[mr]
\bauthor{\bsnm{Erd{\H{o}}s},~\bfnm{P.}\binits{P.}} \AND
\bauthor{\bsnm{R{\'e}nyi},~\bfnm{A.}\binits{A.}}
(\byear{1959}).
\btitle{On random graphs. {I}}.
\bjournal{Publ. Math. Debrecen}
\bvolume{6}
\bpages{290--297}.
\bid{issn={0033-3883}, mr={0120167}}
\end{barticle}
\bptok{imsref}\endbibitem

%b14 ###bbsrt2
\bibitem{ErdosRenyi1960}
\begin{barticle}[mr]
\bauthor{\bsnm{Erd{\H{o}}s},~\bfnm{P.}\binits{P.}} \AND
\bauthor{\bsnm{R{\'e}nyi},~\bfnm{A.}\binits{A.}}
(\byear{1961}).
\btitle{On the evolution of random graphs}.
\bjournal{Bull. Inst. Internat. Statist.}
\bvolume{38}
\bpages{343--347}.
\bid{mr={0148055}}
\end{barticle}
\bptok{imsref}\endbibitem

%b15 ###bbsrt2
\bibitem{EthierKurtz}
\begin{bbook}[mr]
\bauthor{\bsnm{Ethier},~\bfnm{Stewart~N.}\binits{S.N.}} \AND
\bauthor{\bsnm{Kurtz},~\bfnm{Thomas~G.}\binits{T.G.}}
(\byear{1986}).
\btitle{Markov Processes: Characterization and Convergence}.
\bseries{Wiley Series in Probability and Mathematical Statistics:
Probability and Mathematical Statistics}.
\blocation{New York}:
\bpublisher{Wiley}.
\bid{doi={10.1002/9780470316658}, mr={0838085}}
\end{bbook}
\bptok{imsref}\endbibitem

%b16 ###bbsrt2
\bibitem{Faloutsos1999}
\begin{barticle}[auto:STB|2014/02/12|14:17:21]
\bauthor{\bsnm{Faloutsos},~\bfnm{M.}\binits{M.}},
\bauthor{\bsnm{Faloutsos},~\bfnm{P.}\binits{P.}} \AND
\bauthor{\bsnm{Faloutsos},~\bfnm{C.}\binits{C.}}
(\byear{1999}).
\btitle{On power-law relationships of the Internet topology}.
 \bjournal{ACM Comp. Comm. Review}
\bvolume{29}
\bpages{251--262}.
\end{barticle}
\bptok{imsref}\endbibitem

%b17 ###bbsrt2
\bibitem{HannekeFuXing2010}
\begin{barticle}[mr]
\bauthor{\bsnm{Hanneke},~\bfnm{Steve}\binits{S.}},
\bauthor{\bsnm{Fu},~\bfnm{Wenjie}\binits{W.}} \AND
\bauthor{\bsnm{Xing},~\bfnm{Eric~P.}\binits{E.P.}}
(\byear{2010}).
\btitle{Discrete temporal models of social networks}.
\bjournal{Electron. J. Stat.}
\bvolume{4}
\bpages{585--605}.
\bid{doi={10.1214/09-EJS548}, issn={1935-7524}, mr={2660534}}
\end{barticle}
\bptok{imsref}\endbibitem

%b18 ###bbsrt2
\bibitem{Kingman1982}
\begin{barticle}[mr]
\bauthor{\bsnm{Kingman},~\bfnm{J.~F.~C.}\binits{J.F.C.}}
(\byear{1982}).
\btitle{The coalescent}.
\bjournal{Stochastic Process. Appl.}
\bvolume{13}
\bpages{235--248}.
\bid{doi={10.1016/0304-4149(82)90011-4}, issn={0304-4149}, mr={0671034}}
\end{barticle}
\bptok{imsref}\endbibitem

%b19 ###bbsrt2
\bibitem{Kolaczykbook}
\begin{bbook}[mr]
\bauthor{\bsnm{Kolaczyk},~\bfnm{Eric~D.}\binits{E.D.}}
(\byear{2009}).
\btitle{Statistical Analysis of Network Data: Methods and Models}.
\bseries{Springer Series in Statistics}.
\blocation{New York}:
\bpublisher{Springer}.
\bid{doi={10.1007/978-0-387-88146-1}, mr={2724362}}
\end{bbook}
\bptok{imsref}\endbibitem

%b20 ###bbsrt2
\bibitem{LiAldersonDoyleWillingerGraphs}
\begin{barticle}[mr]
\bauthor{\bsnm{Li},~\bfnm{Lun}\binits{L.}},
\bauthor{\bsnm{Alderson},~\bfnm{David}\binits{D.}},
\bauthor{\bsnm{Doyle},~\bfnm{John~C.}\binits{J.C.}} \AND
\bauthor{\bsnm{Willinger},~\bfnm{Walter}\binits{W.}}
(\byear{2005}).
\btitle{Towards a theory of scale-free graphs: Definition, properties,
and implications}.
\bjournal{Internet Math.}
\bvolume{2}
\bpages{431--523}.
\bid{issn={1542-7951}, mr={2241756}}
\end{barticle}
\bptok{imsref}\endbibitem

%b21 ###bbsrt2
\bibitem{LovaszBook}
\begin{bbook}[mr]
\bauthor{\bsnm{Lov{\'a}sz},~\bfnm{L{\'a}szl{\'o}}\binits{L.}}
(\byear{2012}).
\btitle{Large Networks and Graph Limits}.
\bseries{American Mathematical Society Colloquium Publications}
\bvolume{60}.
\blocation{Providence, RI}:
\bpublisher{Amer. Math. Soc.}
\bid{mr={3012035}}
\end{bbook}
\bptok{imsref}\endbibitem

%b22 ###bbsrt2
\bibitem{LovaszSzegedy2006}
\begin{barticle}[mr]
\bauthor{\bsnm{Lov{\'a}sz},~\bfnm{L{\'a}szl{\'o}}\binits{L.}} \AND
\bauthor{\bsnm{Szegedy},~\bfnm{Bal{\'a}zs}\binits{B.}}
(\byear{2006}).
\btitle{Limits of dense graph sequences}.
\bjournal{J. Combin. Theory Ser. B}
\bvolume{96}
\bpages{933--957}.
\bid{doi={10.1016/j.jctb.2006.05.002}, issn={0095-8956}, mr={2274085}}
\end{barticle}
\bptok{imsref}\endbibitem

%b23 ###bbsrt2
\bibitem{McCullagh2002}
\begin{barticle}[mr]
\bauthor{\bsnm{McCullagh},~\bfnm{Peter}\binits{P.}}
(\byear{2002}).
\btitle{What is a statistical model?}
\bjournal{Ann. Statist.}
\bvolume{30}
\bpages{1225--1310}.
\bnote{With comments and a rejoinder by the author}.
\bid{doi={10.1214/aos/1035844977}, issn={0090-5364}, mr={1936320}}
\end{barticle}
\bptok{imsref}\endbibitem

%b24 ###bbsrt2
\bibitem{McCullaghPitmanWinkel2008}
\begin{barticle}[mr]
\bauthor{\bsnm{McCullagh},~\bfnm{Peter}\binits{P.}},
\bauthor{\bsnm{Pitman},~\bfnm{Jim}\binits{J.}} \AND
\bauthor{\bsnm{Winkel},~\bfnm{Matthias}\binits{M.}}
(\byear{2008}).
\btitle{Gibbs fragmentation trees}.
\bjournal{Bernoulli}
\bvolume{14}
\bpages{988--1002}.
\bid{doi={10.3150/08-BEJ134}, issn={1350-7265}, mr={2543583}}
\end{barticle}
\bptok{imsref}\endbibitem

%b25 ###bbsrt2
\bibitem{Milgram1967}
\begin{barticle}[auto:STB|2014/02/12|14:17:21]
\bauthor{\bsnm{Milgram},~\bfnm{S.}\binits{S.}}
(\byear{1967}).
\btitle{The small world problem}.
\bjournal{Psych. Today}
\bvolume{1}
\bpages{60--67}.
\end{barticle}
\bptok{imsref}\endbibitem

%b26 ###bbsrt2
\bibitem{Pitman2005}
\begin{bbook}[mr]
\bauthor{\bsnm{Pitman},~\bfnm{J.}\binits{J.}}
(\byear{2006}).
\btitle{Combinatorial Stochastic Processes}.
\bseries{Lecture Notes in Math.}
\bvolume{1875}.
\blocation{Berlin}:
\bpublisher{Springer}.
\bnote{Lectures from the 32nd Summer School on Probability Theory held
in Saint-Flour, July 7--24, 2002, With a foreword by Jean Picard}.
\bid{mr={2245368}}
\end{bbook}
\bptok{imsref}\endbibitem

%b27 ###bbsrt2
\bibitem{RinaldoShalizi2013}
\begin{barticle}[mr]
\bauthor{\bsnm{Shalizi},~\bfnm{Cosma~Rohilla}\binits{C.R.}} \AND
\bauthor{\bsnm{Rinaldo},~\bfnm{Alessandro}\binits{A.}}
(\byear{2013}).
\btitle{Consistency under sampling of exponential random graph models}.
\bjournal{Ann. Statist.}
\bvolume{41}
\bpages{508--535}.
\bid{doi={10.1214/12-AOS1044}, issn={0090-5364}, mr={3099112}}
\end{barticle}
\bptok{imsref}\endbibitem

%b28 ###bbsrt2
\bibitem{vanderHofstadBook}
\begin{bmisc}[auto:STB|2014/02/12|14:17:21]
\bauthor{\bsnm{van~der Hofstad},~\bfnm{R.}\binits{R.}}
(\byear{2012}).
\bhowpublished{Random walks and complex networks.
Lecture notes, in preparation}.
\end{bmisc}
\bptok{imsref}\endbibitem

%b29 ###bbsrt2
\bibitem{WattsStrogatz1998}
\begin{barticle}[auto:STB|2014/02/12|14:17:21]
\bauthor{\bsnm{Watts},~\bfnm{D.}\binits{D.}} \AND
\bauthor{\bsnm{Strogatz},~\bfnm{S.}\binits{S.}}
(\byear{1998}).
\btitle{Collective dynamics of `small-world' networks}.
\bjournal{Nature}
\bvolume{393}
\bpages{440--442}.
\end{barticle}
\bptok{imsref}\endbibitem

%b30 ###bbsrt2
\bibitem{SixDegrees}
\begin{bbook}[mr]
\bauthor{\bsnm{Watts},~\bfnm{Duncan~J.}\binits{D.J.}}
(\byear{2003}).
\btitle{Six Degrees: The Science of a Connected Age}.
\blocation{New York}:
\bpublisher{W.W. Norton}.
\bid{mr={2041642}}
\bptnote{check year}\end{bbook}
\bptok{imsref}\endbibitem

%b31 ###bbsrt2
\bibitem{WillingerAldersonDoyle2009}
\begin{barticle}[mr]
\bauthor{\bsnm{Willinger},~\bfnm{Walter}\binits{W.}},
\bauthor{\bsnm{Alderson},~\bfnm{David}\binits{D.}} \AND
\bauthor{\bsnm{Doyle},~\bfnm{John~C.}\binits{J.C.}}
(\byear{2009}).
\btitle{Mathematics and the {I}nternet: A source of enormous confusion
and great potential}.
\bjournal{Notices Amer. Math. Soc.}
\bvolume{56}
\bpages{586--599}.
\bid{issn={0002-9920}, mr={2509062}}
\end{barticle}
\bptok{imsref}\endbibitem



\end{thebibliography}
\end{document}